\documentclass[a4paper,oneside,11pt]{article}

\usepackage{amsmath}
\usepackage{amssymb}
\usepackage[utf8x]{inputenc} 
\usepackage{accents}
\usepackage{mathtools}
\usepackage{amsthm}
\usepackage[english]{babel}
\usepackage{letltxmacro}
\usepackage{xparse}
\usepackage[dvipsnames,svgnames,x11names,hyperref]{xcolor}
\usepackage{graphicx}

\usepackage[backgroundcolor = white, linecolor = Cerulean, textcolor = Cerulean, bordercolor = Cerulean]{todonotes}

\usepackage{fullpage}
\usepackage{fancyhdr}
\usepackage{ebgaramond}
\usepackage{euscript}
\usepackage{mathrsfs}
\usepackage{accents}
\usepackage[bbgreekl]{mathbbol}
\usepackage[euler]{textgreek}
\usepackage[scr=boondoxo]{mathalfa}
\DeclareSymbolFontAlphabet{\mathbbm}{bbold}
\DeclareSymbolFontAlphabet{\mathbb}{AMSb}%
\usepackage{mleftright}
\usepackage{dirtytalk}
\usepackage{shuffle}
\usepackage{scalerel}
\usepackage{tensor}
\usepackage{soul}
\setul{0pt}{0.4pt}

\makeatletter
\newcommand*\bigcdot{\mathpalette\bigcdot@{.5}}
\newcommand*\bigcdot@[2]{\mathbin{\vcenter{\hbox{\scalebox{#2}{$\m@th#1\bullet$}}}}}
\makeatother

\mathchardef\mhyphen="2D

\DeclareFontFamily{U}{mathx}{\hyphenchar\font45}
\DeclareFontShape{U}{mathx}{m}{n}{
	<5> <6> <7> <8> <9> <10>
	<10.95> <12> <14.4> <17.28> <20.74> <24.88>
	mathx10
}{}
\DeclareSymbolFont{mathx}{U}{mathx}{m}{n}
\DeclareFontSubstitution{U}{mathx}{m}{n}
\DeclareMathAccent{\widecheck}{0}{mathx}{"71}
\DeclareMathAccent{\wideparen}{0}{mathx}{"75}

\usepackage[ruled]{algorithm2e}

\usepackage{listings}
\DeclarePairedDelimiter{\abs}{\lvert}{\rvert}

\usepackage{tikz}
\usetikzlibrary{decorations.pathreplacing}
\usetikzlibrary{decorations.markings}
\usepackage{pgfplots}
\usepackage{wrapfig}
\usepackage{caption}
\usepackage[percent]{overpic}
\usepackage{tikz-cd}

\usepackage{aliascnt}
\usepackage[pdftex, ocgcolorlinks, colorlinks = true,
linkcolor = imperialBrick,
citecolor = imperialGreen,
urlcolor = imperialProcess
]{hyperref}
\AtBeginDocument{%
	\LetLtxMacro\autoreforig\autoref
	\RenewDocumentCommand{\autoref}{som}{%
		\IfValueT{#2}{[}%
		\IfBooleanTF{#1}{%
			\autoreforig*{#3}%
		}{%
			\autoreforig{#3}%
		}%
		\IfValueT{#2}{,\space#2]}%
	}
}

\usepackage{nameref}

\addto\extrasenglish{
	
} 

\addto\extrasenglish{
	
} 

\addto\extrasenglish{
	
} 

\newcommand{\mynewtheorem}[2]{
	\newaliascnt{#1}{dummy}
	\newtheorem{#1}[#1]{#2}
	\aliascntresetthe{#1}
	\expandafter\def\csname #1autorefname\endcsname{#2}
}

\newtheoremstyle{note}
{\topsep}   
{\topsep}   
{}  
{0pt}       
{\itshape\bfseries} 
{.}         
{5pt plus 1pt minus 1pt} 
{{\color{imperialTangerine}\thmname{#1}\emph{\thmnumber{ #2}}}\hspace{0.1em}\textnormal{\thmnote{ (#3)}}}     

\theoremstyle{plain}
\mynewtheorem{thm}{Theorem}
\mynewtheorem{thmdef}{Theorem/Definition}
\mynewtheorem{prop}{Proposition}
\mynewtheorem{cor}{Corollary}
\mynewtheorem{lem}{Lemma}
\mynewtheorem{conj}{Conjecture}
\theoremstyle{definition}
\mynewtheorem{defn}{Definition}
\mynewtheorem{expl}{Example}
\mynewtheorem{prob}{Problem}
\mynewtheorem{ques}{Question}
\mynewtheorem{cond}{Condition}
\mynewtheorem{conv}{Convention}
\mynewtheorem{defthm}{Definition/Theorem}
\theoremstyle{remark}
\mynewtheorem{rem}{Remark}

\theoremstyle{note}

\definecolor{imperialBlue}{RGB}{0, 62, 116}
\definecolor{imperialBrick}{RGB}{165,25,0}
\definecolor{imperialProcess}{RGB}{0,133,202}
\definecolor{imperialGreen}{RGB}{2,137,59}
\definecolor{imperialRed}{RGB}{221,37,1}
\definecolor{imperialOrange}{RGB}{210,64,0}
\definecolor{imperialBlue2}{RGB}{0,110,175}
\definecolor{imperialTangerine}{RGB}{236,115,0}
\definecolor{imperialPurple}{RGB}{101,48,152}
\definecolor{imperialLime}{RGB}{196,214,0}
\definecolor{imperialKermit}{RGB}{102,164,10}

\definecolor{kermit}{RGB}{102,164,10}
\definecolor{teal}{RGB}{0,142,170}
\definecolor{tangerine}{RGB}{236,115,0}
\definecolor{raspberry}{RGB}{145,0,72}
\definecolor{lime}{RGB}{196,214,0}

\colorlet{chaptergrey}{imperialBlue}

\usepackage{cutwin,kantlipsum,mwe}
\usepackage{regexpatch}
\opencutright

\makeatletter
\xpretocmd{\cutout}{\leavevmode\hrule \@height\z@ \@width\linewidth\relax}
\makeatother

\usepackage{iftex}
\iftutex
\usepackage{unicode-math}
\setmathfontface\altgrfont{GFS Artemisia Italic}[Scale=MatchLowercase]
\newcommand{\text{\textdelta}}{\altgrfont{δ}}
\else
\usepackage[T1]{fontenc}
\DeclareSymbolFont{altgr}{OML}{antt}{m}{it}
\DeclareMathSymbol{\sko }{\mathord}{altgr}{"0E}
\fi

\let\oldKronecker\sko
\renewcommand{\sko}{\oldKronecker \hspace{-0.1em}}

\def\kron{\text{\textdelta}}

\makeatletter
\def\upintkern@{\mkern-7mu\mathchoice{\mkern-3.5mu}{}{}{}}
\def\upintdots@{\mathchoice{\mkern-4mu\@cdots\mkern-4mu}%
	{{\cdotp}\mkern1.5mu{\cdotp}\mkern1.5mu{\cdotp}}%
	{{\cdotp}\mkern1mu{\cdotp}\mkern1mu{\cdotp}}%
	{{\cdotp}\mkern1mu{\cdotp}\mkern1mu{\cdotp}}}

\newcommand{\UpMultiIntegral}[1]{%
	\edef\ints@c{\noexpand\upintop
		\ifnum#1=\z@\noexpand\upintdots@\else\noexpand\upintkern@\fi
		\ifnum#1>\tw@\noexpand\upintop\noexpand\upintkern@\fi
		\ifnum#1>\thr@@\noexpand\upintop\noexpand\upintkern@\fi
		\noexpand\upintop
		\noexpand\ilimits@
	}%
	\futurelet\@let@token\ints@a
}
\makeatother

\DeclareFontFamily{OMX}{mdbch}{}
\DeclareFontShape{OMX}{mdbch}{m}{n}{ <->s * [0.8]  mdbchr7v }{}
\DeclareFontShape{OMX}{mdbch}{b}{n}{ <->s * [0.8]  mdbchb7v }{}
\DeclareFontShape{OMX}{mdbch}{bx}{n}{<->ssub * mdbch/b/n}{}

\DeclareSymbolFont{uplargesymbols}{OMX}{mdbch}{m}{n}
\SetSymbolFont{uplargesymbols}{bold}{OMX}{mdbch}{b}{n}
\DeclareMathSymbol{\upintop}{\mathop}{uplargesymbols}{82}
\DeclareMathSymbol{\upointop}{\mathop}{uplargesymbols}{"48}

\DeclareFontEncoding{MDB}{}{}
\DeclareFontFamily{MDB}{mdbch}{}
\DeclareFontShape{MDB}{mdbch}{m}{n}{ <->s * [0.8]  mdbchrmb }{}
\DeclareFontShape{MDB}{mdbch}{b}{n}{ <->s * [0.8]  mdbchbmb }{}
\DeclareFontShape{MDB}{mdbch}{bx}{n}{<->ssub * mdbch/b/n}{}
\DeclareFontSubstitution{MDB}{cmr}{m}{n}
\DeclareSymbolFont{mathdesignB}{MDB}{mdbch}{m}{n}%
\SetSymbolFont{mathdesignB}{bold}{MDB}{mdbch}{b}{n}%
\DeclareMathSymbol{\upintclockwise}{\mathop}{mathdesignB}{128}
\DeclareMathSymbol{\upointclockwise}{\mathop}{mathdesignB}{130}
\DeclareMathSymbol{\upointctrclockwise}{\mathop}{mathdesignB}{132}
\DeclareMathSymbol{\upoiint}{\mathop}{mathdesignB}{134}
\DeclareMathSymbol{\upoiiint}{\mathop}{mathdesignB}{136}

\makeatletter
\newcommand{\upint}{\DOTSI\upintop\ilimits@}
\newcommand{\upoint}{\DOTSI\upointop\ilimits@}
\makeatother

\renewcommand{\int}{\upint}

\usepackage[scr=boondoxo]{mathalfa}

\newcommand{\ddif}{\mathbbm d}

\newcommand{\bbT}{\mathbb T}
\newcommand{\bbF}{\mathbb F}

\newcommand{\bbR}{\mathbb R}

\usepackage[euler]{textgreek}

\newcommand{\dif}{\mathrm{d}}

\title{\textsc{Projections of SDEs onto Submanifolds}}

\author{John Armstrong \\ King's College London \\ \href{mailto:john.armstrong@kcl.ac.uk}{\nolinkurl{john.armstrong@kcl.ac.uk}} \and Damiano Brigo \\ Imperial College London \\ \href{mailto:damiano.brigo@imperial.ac.uk}{\nolinkurl{damiano.brigo@imperial.ac.uk}} \and Emilio Ferrucci \\ University of Oxford \\ \href{mailto:emilio.rossiferrucci@maths.ox.ac.uk}{\nolinkurl{emilio.rossiferrucci@maths.ox.ac.uk}}}

\date{\today}

\makeatletter
\newcommand*{\starsection}[1]{%
	\section*{#1}%
	\NR@gettitle{#1}%
}
\makeatother

\begin{document}

\maketitle

\abstract{In \cite{ABR19} the authors define three projections of $\bbR^d$-valued stochastic differential equations (SDEs) onto submanifolds: the Stratonovich, It\^o-vector and It\^o-jet projections. In this paper, after a brief survey of SDEs on manifolds, we begin by giving these projections a natural, coordinate-free description, each in terms of a specific representation of manifold-valued SDEs. We proceed by deriving formulae for the three projections in ambient $\bbR^d$-coordinates. We use these to show that the It\^o-vector and It\^o-jet projections satisfy respectively a weak and mean-square optimality criterion \say{for small t}: this is achieved by solving constrained optimisation problems. These results confirm, but do not rely on the approach taken in \cite{ABR19}, which is formulated in terms of weak and strong It\^o-Taylor expansions. In the final section we exhibit examples showing how the three projections can differ, and explore alternative notions of optimality.}

\section*{Introduction} \label{sec:projIntro}
\addcontentsline{toc}{section}{Introduction}

Consider the following problem: we are given an autonomous ODE
\begin{equation}\label{ODE}
	\dot X_t = F(X_t)
\end{equation}
in $\bbR^d$, and a smooth embedded manifold $M \hookrightarrow \bbR^d$. Let $\pi$ be the metric projection of a tubular neighbourhood of $M$ onto $M$ (see \eqref{defpi} below). We seek an $M$-valued ODE, i.e.\ a vector field $\overline F$ on $M$, tangent at each point to $M$, with the property that the solution to 
\begin{equation}
	\dot Y_t = \overline F(Y_t)
\end{equation}
is optimal in the sense that the first coefficient of the Taylor expansion in $t = 0$ of either
\begin{equation}\label{ODEproj}
	|Y_t - X_t|^2 \quad \text{or} \quad |Y_t - \pi(X_t)|^2
\end{equation}
is minimised for any initial condition $X_0 = Y_0 = y_0 \in M$. This requirement represents the slowest possible divergence of $Y$ from the original solution $X$ (resp. from its metric projection on $M$), subject to the constraint of $Y$ arising as the solution of a closed form ODE on $M$. It is an easy exercise (using \eqref{JPQ} below) to check that these optimisation problems both result in the same solution, which consists in $\overline F(y)$ being the orthogonal projection of the vector $F(y)$ onto the tangent space $T_yM$.

The paper \cite{ABR19}, which is motivated by applications to non-linear filtering, explores an extension of this problem to the case of SDEs. The optimality criteria \eqref{ODEproj} do not carry over in a straightforward fashion, and are formulated through the machinery of weak and strong It\^o-Taylor expansions. In this chapter we tackle the same problem through a different perspective, which we proceed to describe.

In \autoref{sec:SDEs} we begin with a survey of SDEs on manifolds. Here we introduce three ways of representing them: the Stratonovich, Schwartz-Meyer (or 2-jet) and It\^o representations. The first and second have the advantage of not requiring a connection on the tangent bundle of the manifold, the second and third are defined in terms of the It\^o integral, while the first and third have vector coefficients. Focusing on the diffusion case, we show how to pass from one representation to another. In \autoref{sec:projEmbedded} we prepare the framework for manifolds $M$ embedded in $\bbR^d$. These are entirely general Riemannian manifolds, due to the Nash embedding theorem, and have the advantage of being describable using ambient coordinates. We use this framework to study the equations introduced in the previous section, on embedded manifolds. In \autoref{sec:projecting} we associate to each manifold-valued SDE representation a natural projection, which gives rise to an SDE on a submanifold: the Stratonovich projection (defined by projecting the Stratonovich coefficients), the It\^o-jet projection (defined by projecting the Schwartz morphism, or 2-jet, which defines the SDE), and the It\^o-vector projection (defined by projecting the It\^o coefficients, and interpreting the resulting equation w.r.t.\ the Riemannian connection on the embedded submanifold). These projections coincide with the ones introduced in \cite{ABR19}, but are given a more solid theoretical underpinning, which sheds light on their analytic and probabilistic properties. We then derive formulae for the three projections, preferring ambient coordinates to local coordinates. In \autoref{sec:optimal} we formulate the optimality criteria satisfied by the It\^o-vector and It\^o-jet projections using respectively an explicit weak and mean-square formulation, instead of invoking It\^o-Taylor expansions as done in \cite{ABR19}. This has the advantage of representing a more tangible property of the solution, and is accompanied by an argument, based on martingale estimates, used to deal with the problem of the solution exiting the tubular neighbourhood of $M$. Our main theorems \autoref{maintheoremV} and \autoref{maintheorem} replicate the findings \cite[Theorem 4.4 and Theorem 4.7]{ABR19} in this new setting. The fact that the Stratonovich projection does not satisfy either of these optimality criteria is a confirmation of the fact that It\^o calculus on manifolds can be of great interest. In \autoref{sec:further} we provide examples showing that the three projections are genuinely distinct, we prove the It\^o projections are optimal also when formulating the optimality criteria using $M$'s intrinsic geometry, and explore notions of optimality that are satisfied by the more na\"ive Stratonovich projection.

Although the material presented here overlaps to a significant degree with the ideas of \cite{ABR19}, this paper --- the contents of which also appear in the third author's PhD thesis \cite{Fer22} --- is entirely self-contained. Moreover, we believe the framework chosen here has a number of advantages of which we hope to make use in future work, as described in \nameref{sec:projConcl}. 

\section{SDEs on manifolds} \label{sec:SDEs}
We begin this chapter with a primer on manifold-valued SDEs. Since manifolds, unlike Euclidean space, do not come naturally equipped with coordinates, especially not global ones, the challenge is to express an SDE using intrinsic, coordinate-free notions. Equivalently, one can define an SDE locally in an arbitrary chart, and show that the property of a process of being a solution does not depend on the chart. The coordinate-free definition of a time-homogeneous ODE on a smooth, $m$-dimensional manifold $M$ is well known: this consists of a tangent vector field, i.e.\ a section of the tangent bundle of $M$, $V \in \Gamma TM$. We will denote $\Gamma$ the set of sections of a fibre bundle, i.e.\ the smooth right inverses to the bundle projection. A solution to the ODE defined by $V$ is a smooth curve $X$, defined on some interval of $\bbR$, with the property that $\dot X_t = V_{X_t}$ for all $t$. This is a coordinate-free definition, and in a chart $\varphi \colon U \to \bbR^m$ ($U$ open set in $M$) it corresponds to requiring that, writing $\varphi(X_t) = {^\varphi \!}X_t$ and $V_x = {^\varphi }V_x^k \partial_x \varphi_k$, we have ${^\varphi \!}\dot X_t^k = {^\varphi}V^k_{X_t}$ for all $t$ for which both sides are defined. Notice the sum over $k$: this is the Einstein summation convention, which we will use throughout this thesis whenever possible; also, $\partial_x\varphi_k$ are the elements of the basis of $T_xM$ defined by the chart $\varphi$:
\begin{equation}
	\partial_x\varphi_k(f) \coloneqq \frac{\partial (f \circ \varphi^{-1})}{\partial x^k}(\varphi(x)) \quad \text{for } f \in C^\infty M
\end{equation}
In this section we will give similar descriptions of Stratonovich and It\^o (non path-dependent) SDEs on manifolds. From now on we will avoid the $\varphi$ superscripts when no ambiguity occurs, e.g.\ the previous identity will be written $\dot X_t^k = V^k_{X_t}$.

We begin with the Stratonovich case, following mainly \cite[Ch.\ VII]{E89}, although the topic is well known. As for the familiar $\bbR^d$-valued case we 
will also need a driving semimartingale, which, given the context we are working in can be taken to be valued in another manifold $N$, of dimension $n$. Given a stochastic setup $(\Omega, \mathcal F_\cdot, P)$ satisfying the usual conditions, a continuous adapted stochastic process $Z \colon \Omega \times \bbR_{\geq 0} \to N$ is said to be a \emph{semimartingale} if, for all $f \in C^\infty N$, $f(Z)$ is a semimartingale. Just as for the ODE case, what is needed to define a Stratonovich SDE in $M$ driven by $Z$ is a section of some vector bundle: in this case, however, the bundle is no longer just $TM$, but $\mathrm{Hom}(TN,TM) \to M \times N$, i.e.\ the vector bundle of linear maps from $TN$ to $TM$. An element $F \in \Gamma \mathrm{Hom}(TN,TM)$ corresponds to a smooth map $M \times N \ni (x,z) \mapsto F(x,z) \in \mathrm{Hom}(T_zN, T_xM)$. The Stratonovich SDE
\begin{equation}\label{StratonovichSDE}
	\dif X_t = F(X_t,Z_t) \circ \dif Z_t
\end{equation}
in local coordinates (this requires choosing a chart both on $N$ and on $M$) as $\dif X_t^k = F^k_\gamma(X_t, Z_t) \circ \dif Z^\gamma_t$ on random intervals that make both sides of the expression well defined. We will always use Greek letters as indices for the driving process, and Latin letters as indices for the solution. The key property that allows one to prove that the coordinate formulation of Stratonovich SDEs holds for all other charts (on the intersection of their respective domains) is that Stratonovich equations satisfy the first order chain rule: clearly \eqref{StratonovichSDE} would not be similarly well defined with It\^o integration.
One can also define a solution without invoking charts: this entails defining a Stratonovich integral taking as integrator an $M$-valued semimartingale $X$ and as integrand a previsible process $H$ with values in the cotangent bundle of $M$ and relatively compact image (\emph{locally bounded}), s.t.\ at each $t$, $H_t$ is in the fibre at $X_t$: this yields an $\bbR$-valued semimartingale which we can write as
\begin{equation}
	\int_0^\cdot \langle H_s, \circ \dif  X_s \rangle 
\end{equation}
The angle brackets refer to dual pairing of vectors and covectors. This integral is characterised as being the unique map satisfying the following three properties
\begin{description}
	\item[Additivity.] For all locally bounded previsible $H,G$ above $X$ $$ \int_0^\cdot \langle H_s + G_s, \circ \dif X_s \rangle = \int_0^\cdot \langle  H_s, \circ \dif X_s \rangle + \int_0^\cdot \langle  G_s, \circ \dif X_s \rangle$$
	\item[Associativity.] For a real-valued, locally bounded adapted process $\lambda$ $$\int_0^\cdot \langle \lambda_s H_s, \circ \dif X_s \rangle  =  \int_0^\cdot \lambda_s \circ \dif \textstyle \int_0^s \langle  H_u, \circ \dif X_u \rangle$$
	\item[Change of variable formula.] For all $f \in C^\infty M$ $$\int_0^\cdot \langle \dif_{X_s}f, \circ\dif X_s \rangle = f(X) - f(X_0)$$
\end{description}
where $\dif f$ is the one-form given by taking the differential of $f$. One can then use this integral to say that $X$ solves \eqref{StratonovichSDE} if for all admissible integrands $H$ (even just those arising as the evaluation of a one-form at $X$) 
\begin{equation}\label{fromintegraltoSDE}
	\int_0^\cdot \langle H_s, \circ \dif X_s \rangle = \int_0^\cdot \langle F(X_s,Z_s)^* H_s,  \circ \dif Z_s \rangle 
\end{equation}
where the $*$ denotes dualisation.
\begin{rem}[Autonomousness and explicitness]
	If $N = \bbR^n$ we can call \eqref{StratonovichSDE} \emph{autonomous} if $F(z,x)$ does not depend on $z$, and if $M = \bbR^m$ we can call it 		\emph{explicit} if $F(z,x)$ does not depend on $x$. However, in the general manifold setting these two concepts do not carry over, at least not unless $N$ (resp. $M$) is parallelisable, with a chosen trivialisation of its tangent bundle. An analogous consideration applies to other flavours of SDEs introduced in this section.
\end{rem}
\begin{expl}[Stratonovich diffusion]\label{stratdiffusion}
	An important example is the case where $N = \bbR_{\geq 0} \times \bbR^n$ and $Z_t = (t,W_t)$, $W$ an $n$-dimensional Brownian motion, and $F$ not depending explicitly on $W$. This means \eqref{StratonovichSDE} becomes 
	\begin{equation}\label{stratonovichequation}
		\dif X_t = \sigma_\gamma(X_t, t) \circ \dif W^\gamma_t + b(X_t,t) \dif t
	\end{equation}
	for $\sigma_\gamma, b \in \Gamma \mathrm{Hom}(T\bbR_{\geq 0},TM) = C^\infty(\bbR_{\geq 0}, \Gamma TM)$, $\gamma = 1,\ldots, n$. Stratonovich diffusions are sections of the vector bundle
	\begin{align}\label{stratdiffbundle}
		\begin{split}
			\mathrm{Diff}^{\, n}_\mathrm{Strat} &\coloneqq \{F \in \mathrm{Hom}(T(\bbR_{\geq 0} \oplus \bbR^n), TM) : \forall w_1,w_2 \in \bbR^n \ F(t,w_1;x) = F(t,w_2;x) \} \\
			&\to M \times \bbR_{\geq 0}
		\end{split}
	\end{align}
	i.e.\ elements of the vector space $\Gamma \mathrm{Diff}^{\, n}_\mathrm{Strat}$. Notice that the base space is not $M \times (\bbR_{\geq 0} \times \bbR^n)$, since independence of the Brownian motion allows us to forget the $\bbR^n$ component.
	
\end{expl}
We note that no additional structure on $N$ and $M$, apart from their smooth atlas, is needed to define Stratonovich equations. Stratonovich SDEs are the most used in stochastic differential geometry, as they behave well w.r.t.\ notions of first order calculus: for instance, if there exists an embedded submanifold $M'$ of $M$ such that $F(y,z)$ maps to $T_y M'$ for all $z \in N$ and all $y \in M'$, then the solution to the Stratonovich SDE defined by $F$ started on $M'$ will remain on $M'$ for the duration of its lifetime. This is evident from our intrinsic approach, by considering $F|_{M' \times N}$, but some authors who develop Stratonovich calculus on manifolds extrinsically prove this by showing that the distance between the solution and the manifold (embedded in Euclidean space) is zero \cite[Prop. 1.2.8]{Hsu02}. The existence and uniqueness of solutions to Stratonovich SDEs can be treated by using the Whitney embedding theorem to embed $N$ and $M$ in Euclidean spaces of high enough dimension, and smoothly extending $F$ so that it vanishes outside a compact set containing the manifolds. Invoking the usual existence and uniqueness theorem (e.g.\ \cite[Theorems 38-40]{P05}), and the good behaviour of Stratonovich SDEs w.r.t.\ submanifolds, immediately proves that a unique solution exists up to a positive stopping time, provided $F$ is smooth. We will mostly not be concerned with global-in-time existence in this thesis, although sufficient conditions for such behaviour can usually be obtained by requiring global Lipschitz continuity w.r.t.\ complete Riemannian metrics.

We now pass to It\^o theory on manifolds, as developed in \cite[Ch.VI]{E89}. The difficulty lies in the second order chain rule of the It\^o integral. For this reason, we need to invoke structures of order higher than 1. Let the \emph{second order tangent bundle} of $M$, $\bbT M$, denote the bundle of second order differential operators without a constant term, i.e.\ given a local chart $\varphi$ containing $x$ in its domain, an element of $L_x \in \bbT_x M$ consists of a map
\begin{equation}
	L_x \colon C^\infty M \to \bbR, \quad L_x f = L^k_x \frac{\partial f}{\partial \varphi^k} + L^{ij}_x \frac{\partial^2 f}{\partial \varphi^i \partial \varphi^j}
\end{equation}
The coefficients $L^k_x$, $L^{ij}_x$ obviously depend on $\varphi$, but their existence does not; moreover, requiring $L^{ij}_x = L^{ji}_x$ ensures their uniqueness for the given chart $\varphi$. Note that if the $L_x^{ij}$'s vanish $L_x \in T_xM$. $\bbT M$ is given the unique topology and smooth structure that makes the projection $\bbT M \to M$, $L_x \mapsto x$ a locally trivial surjective submersion. Just as for the first order case, there is an obvious notion of induced bundle map $\bbT f \colon \bbT N \to \bbT M$ for $f \in C^\infty(N,M)$. A chart $\varphi$ containing $x$ in its domain defines the basis 
\begin{equation}
	\{\partial_x \varphi_k, \partial^2_x \varphi_{ij} = \partial^2_x \varphi_{ji} \mid k,i,j = 1,\ldots , n\}
\end{equation}
so the dimension of $\bbT M$ (as a vector bundle) is $m + m(m+1)/2$.
The fundamental properties of $\bbT M$ are summarised the short exact sequence of vector bundles over $M$
\begin{equation}\label{SES}
	\begin{tikzcd}[row sep = tiny]
		0 \arrow[r] &TM \arrow[r,"\mathscr i"] &\bbT M \arrow[r,"\mathscr p"] &TM \odot TM \arrow[r] & 0
	\end{tikzcd}
\end{equation}
with the third term denoting symmetric tensor product, the first map the obvious inclusion and the second map given by
\begin{equation}
	L_x \mapsto \big(f,g \mapsto \frac 12 (L_x(fg) - f(x)L_xg - g(x)L_xf)\big)	
\end{equation}
Roughly speaking, this means that $\bbT M$ is \say{noncanonically the direct sum of $TM$ and $TM \odot TM$}. This short exact sequence of course dualises to a short exact sequence of dual bundles. Elements of $\bbT^*_xM$ can always be represented as $\ddif_x f$, defined by
\begin{equation}\label{ddiff}
	\langle \ddif_x f, L_x \rangle \coloneqq L_x(f)
\end{equation}
for some $f \in C^\infty M$ (this is of course only true at a point: not all sections of $\bbT M$ are of the form $\ddif f$). We now wish to define an It\^o-type equation using second order tangent bundles instead of ordinary tangent bundles. For this we need a notion of field of maps $\bbF (x,z) \colon \bbT_z N \to \bbT_x M$. Since the bundles in question are linear, it is tempting to allow $\bbF(x,z)$ to be an arbitrary linear map, but a more stringent condition is necessary to guarantee well-posedness: the correct requirement is that $\bbF(x,z)$ define a morphism of short exact sequences, i.e.\ a commutative diagram
\begin{equation}\label{Schcondition}
	\begin{tikzcd}
		0 \arrow[r] &T_zN \arrow[d, "F(x{,}z)"] \arrow[r] &\bbT_z N \arrow[d, "\bbF (x{,}z)"] \arrow[r] &T_zN \odot T_zN \arrow[d, "F(x{,}z) \otimes F(x{,}z)"] \arrow[r] & 0 \\
		0 \arrow[r] &T_xM \arrow[r] &\bbT_x M \arrow[r] &T_xM \odot T_xM \arrow[r] & 0
	\end{tikzcd}
\end{equation}
with $F(x,z) = \bbF(x,z)|_{T_zN}$. $\bbF(x,z)$ is then called a \emph{Schwartz morphism}, and we can then view $\bbF$ as being the section of a sub-fibre bundle $\mathrm{Sch}(N,M)$ of $\mathrm{Hom}(\bbT N, \bbT M)$ over $M \times N$ consisting of such maps, which we call the \emph{Schwartz bundle}. Note that $\mathrm{Sch}(N,M)$ is not closed under sum and scalar multiplication taken in the vector bundle $\mathrm{Hom}(\bbT N, \bbT M)$, and thus can only be treated as a fibre bundle. Now, given $\bbF \in \Gamma \mathrm{Sch}(N,M)$, we will give a meaning to the SDE
\begin{equation}\label{SMequation}
	\ddif X_t = \bbF(X_t, Z_t) \ddif Z_t
\end{equation}
which we will call a \emph{Schwartz-Meyer equation}. Heuristically, if $X$ is an $M$-valued semimartingale the second order differential $\ddif X_t$ should be interpreted in local coordinates $\varphi$ as 
\begin{equation}\label{ddif}
	\ddif X_t = \dif X_t^k \partial_{X_t}\varphi_k + \tfrac 12 \dif[X^i,X^j]_t \partial_{X_t}^2 \varphi_{ij} \in \bbT_{X_t}M
\end{equation}
where the first differential is an It\^o differential; this expression is seen to be invariant under change of charts, thanks to the It\^o formula. Then, given charts $\varphi$ in $M$ and $\vartheta$ on $N$, and writing 
\begin{align}
	\begin{split}
		\bbF(x,z)\partial_z\vartheta_\gamma &= \bbF_\gamma^k(x,z) \partial_x \varphi_k + \bbF_\gamma^{ij}(x,z) \partial_x^2 \varphi_{ij} \\
		\bbF(x,z)\partial_z^2 \vartheta_{\alpha \beta} &= \bbF_{\alpha\beta}^k(x,z) \partial_x \varphi_k + \bbF_{\alpha\beta}^{ij}(x,z) \partial_x^2 \varphi_{ij}
	\end{split}
\end{align}
\eqref{SMequation} becomes the system
\begin{equation}\label{SMcases}
	\begin{cases}
		\dif X^k_t = \bbF^k_\gamma (X_t,Z_t) \dif Z_t^\gamma + \frac 12 \bbF^k_{\alpha\beta}(X_t,Z_t) \dif [Z^\alpha, Z^\beta]_t \\
		\frac 12 \dif [X^i,X^j]_t = \bbF^{ij}_\gamma(X_t,Z_t) \dif Z^\gamma_t + \frac 12  \bbF^{ij}_{\alpha\beta}(X_t,Z_t) \dif[Z^\alpha,Z^\beta]_t
	\end{cases}
\end{equation}
Computing the quadratic covariation matrix of $X$ from the first equation above, using the Kunita-Watanabe identity, and comparing with the second results in the requirement that
\begin{equation}\label{Schconditioncoo}
	\bbF^{ij}_\gamma \equiv 0; \quad \bbF_{\alpha\beta}^{ij} \equiv \tfrac 12 \big(\bbF_\alpha^i   \bbF_\beta^j + \bbF_\alpha^j  \bbF_\beta^i \big)
\end{equation}
which correspond precisely to the Schwartz condition \eqref{Schcondition}, and justifies this requirement. \eqref{SMcases} now reduces to its first line, i.e.\ the It\^o SDE
\begin{equation}\label{SMincoords}
	\dif X_t^k = \bbF^k_\gamma (X_t,Z_t) \dif Z_t^\gamma + \tfrac 12 \bbF_{\alpha\beta}^k(X_t,Z_t) \dif[Z^\alpha, Z^\beta]_t
\end{equation}
on random intervals that make both sides of the expression well-defined.
\begin{expl}[Schwartz-Meyer diffusion]\label{SchMeydiffusion}
	Proceeding as in \autoref{stratdiffusion}, but with Schwartz-Meyer equations, we can define the Schwartz-Meyer SDE
	\begin{align}\label{schwartzdiffusion}
		\begin{split}
			\ddif X_t &= \bbF(X_t,t) \ddif Z_t \\
			&= \sigma_\gamma(X_t,t) \dif W^\gamma_t + \bigg(\bbF_0 + \frac 12 \sum_{\gamma = 1}^n \bbF_{\gamma\gamma} \bigg)(X_t,t) \dif t
		\end{split}
	\end{align}
	where we can call $\bbF_\gamma = \sigma_\gamma$ the diffusion coefficients, since they are elements of $C^{\infty}(\bbR_{\geq 0},\Gamma TM)$; this also holds for $\gamma = 0$, but not for $\bbF_{\alpha\beta} \in C^{\infty}(\bbR_{\geq 0},\Gamma \bbT M)$. Therefore the coefficient of $\dif t$, the \say{drift}, cannot be interpreted as a vector. Note that setting $\bbF_{\gamma\gamma} \equiv 0$ does not guarantee that such coefficients will vanish w.r.t.\ another chart, since the transformation rule for them involves the $\bbF_{\alpha\beta}^{ij}$'s which cannot vanish by the second Schwartz condition \eqref{Schconditioncoo}; in other words, there is no way to do away with the non vector-valued drift in \eqref{schwartzdiffusion}. We can consider Schwartz Meyer diffusions as being sections of the fibre bundle
	\begin{align}\label{schdiffbundle}
		\begin{split}
			\mathrm{Diff}^{\, n}_\mathrm{Sch} M &\coloneqq \frac{\{\bbF \in \mathrm{Sch}(\bbR_{\geq 0} \times \bbR^n, M) : \forall w_1,w_2 \in \bbR^n \ \bbF(t,w_1;x) = \bbF(t,w_2;x) \}}{\bbF \sim \mathbb G \Leftrightarrow \bbF_{\gamma\geq 1} = \mathbb G_\gamma, \ \bbF_0 + \frac 12 \sum_{\gamma = 1}^n \bbF_{\gamma\gamma} = \mathbb G_0 + \frac 12 \sum_{\gamma = 1}^n \mathbb G_{\gamma\gamma}} \\ 
			&\to M \times \bbR_{\geq 0}
		\end{split}	
	\end{align}
	This means that, similarly to the case of \eqref{stratdiffbundle} we are only considering $\bbF$'s that do not depend explicitly on the Brownian motion, and we are quotienting out the part that is not relevant for \eqref{schwartzdiffusion}.
\end{expl}

Just as for Stratonovich SDEs, Schwartz-Meyer equations can also be seen to come from an integral
\begin{equation}\label{SMintegral}
	\int_0^\cdot \langle \mathbb H_s, \ddif X_s \rangle	
\end{equation}
where the process $\mathbb H$ is now valued in $\bbT^* M$. The axioms for this \emph{Schwartz-Meyer integral} are similar:
\begin{description}
	\item[Additivity.] For all locally bounded previsible $\mathbb H, \mathbb G$ above $X$ $$ \int_0^\cdot \langle \mathbb  H_s + \mathbb G_s,  \ddif X_s \rangle = \int_0^\cdot \langle  \mathbb H_s,  \ddif X_s \rangle + \int_0^\cdot \langle \mathbb G_s,  \ddif X_s \rangle$$
	\item[Associativity.] For a real-valued, locally bounded adapted process $\lambda$ $$\int_0^\cdot \langle \lambda_s \mathbb H_s,  \ddif X_s \rangle  =  \int_0^\cdot \lambda_s \dif \textstyle \int_0^s \langle  H_u, \ddif X_u \rangle$$
	\item[Change of variable formula.] For all $f \in C^\infty M$ $$\int_0^\cdot \langle \ddif_{X_s}f, \ddif X_s \rangle = f(X) - f(X_0)$$
\end{description}
Notice how It\^o integration is used in the associativity axiom. The property of a process of being a solution of \eqref{SMequation} is then defined in complete analogy to \eqref{fromintegraltoSDE}.

The recent paper \cite{AB18} treats SDEs on manifolds using a representation which is similar to that of \eqref{SMequation}, but which has a distinct advantage when it comes to numerical schemes. Here the authors focus on the autonomous diffusion case, without explicitly taking time as a driver ($N = \bbR^n$, $Z_t = W_t$), and take the field of Schwartz morphisms $\bbF$ to be \emph{induced} by a \emph{field of maps} i.e.\ a smooth function $f \colon \bbR^n \times M \to M$, $f_x \coloneqq f(\cdot,x)$, s.t. for all $x \in M$, $f_x(0) = x$: this means 
\begin{equation}\label{jets}
	\bbF(x) = \bbT_0 f_x
\end{equation}
In coordinates $\varphi$ on $M$ this amounts to
\begin{equation}
	\sigma^k_\gamma(x) = \frac{\partial (\varphi^k \circ f_x)}{\partial w^\gamma}(0), \quad \bbF^k_{\alpha\beta}(x) = \frac{\partial^2 (\varphi^k \circ f_x)}{\partial w^\alpha \partial w^\beta} (0)
\end{equation}
with $\bbF_0 = 0$ (note how the drift comes from the quadratic variation of Brownian motion, without having to require time as a driving process). This particular form of $\bbF$ is useful because it automatically defines a numerical scheme for the solution of the SDE, similar to the Euler scheme, which cannot be defined in a coordinate-free way on a manifold: the linear structure lacked by $M$ is replaced with iterative interpolations along the $f_x$'s. This also has the advantage of guaranteeing that if the maps are valued in $M$, so are all the approximations.

\say{It\^o-type} Diffusions on manifolds have also been investigated by other authors, most notably by \cite[Ch.4]{BeDa90} (although we refer to the more recent exposition \cite[\S 7.2]{Gl11}), who call the bundle $\mathrm{Diff}_\mathrm{Sch}^{\, n} M$ the \emph{It\^o bundle}, and give a local description of it. Although we will not need this formulation in the following sections, we include a description of it to establish the link with the other approach. There are (at least) two ways of describing a fibre bundle $\pi \colon E \to M$: one is by simply exhibiting the manifolds $E,M$ and the surjective submersion $\pi$, and by checking local triviality; this is the approach taken here. The second approach involves declaring the base space $M$, the structure group $G$ (a Lie group), the typical fibre $F$ (a smooth manifold, carrying a left action of $G$ by smooth maps) and a covering $\{U_\lambda\}_\lambda$ of $M$ together with maps $g_{\nu\mu} \colon U_\mu \cap U_\nu \to G$ satisfying the cocycle conditions $\forall \lambda, \mu, \nu \ g_{\nu\mu}g_{\mu\lambda} = g_{\nu\lambda}$. Then the total space and bundle projection can be reconstructed by gluing all the $U_\lambda \times F$'s together according to the $g_{\nu\mu}$'s:
\begin{equation}\label{gluing}
	E \coloneqq \frac{\bigcup_{\lambda} \{\lambda \} \times U_\lambda \times F}{(\mu,x,e) \sim (\nu,y,f) \Leftrightarrow x = y, \ f = g_{\nu\mu}(x).e} \xrightarrow{\pi} M, \quad [\mu, x, e] \mapsto x
\end{equation}
Of course, the local description can be obtained from the ordinary one by fixing a local trivialisation, a model for the fibre, a Lie group capturing all transformations of the fibres, etc. Now, we define the candidate bundle of Schwartz-Meyer diffusions to have base space $M \times \bbR_{\geq 0}$ and typical fibre $\mathrm{Hom}(\bbR^n,\bbR^m) \oplus \bbR^m$. Recall that we observed that the Schwartz bundle is not linear: this should rule out the usual choices $G = GL(n,\bbR), O(n)$, valid for vector bundles. Indeed, the transformation laws for $\mathrm{Diff}_\mathrm{Sch}^{\, n} M$ are succinctly modelled by the \emph{It\^o group}
\begin{align}
	\mathfrak I^m &\coloneqq GL(m,\bbR) \times \mathrm{Hom}(\bbR^m \odot \bbR^m, \bbR^m) \\
	(A,a)(B,b) &\coloneqq (A \circ B, A \circ b + a \circ (B \otimes B))
\end{align}
with identity $(I_m,0)$, acting on $\mathrm{Hom}(\bbR^n,\bbR^m) \oplus \bbR^m$ from the left by
\begin{equation}
	(A,a).(\sigma,\eta) \coloneqq \big( A \circ \sigma, A \eta + \tfrac 12 \mathrm{tr}( a \circ (\sigma \otimes \sigma) ) \big)
\end{equation}
where the trace is taken componentwise. Given an open covering $\{U_\lambda\}_\lambda$ (consisting of, say, open balls) of $M$, and charts $\varphi_\lambda \colon U_\lambda \to \bbR^m$, we define
\begin{equation}
	g_{\nu\mu}(x \in U_\mu \cap U_\nu) \coloneqq \big( J(\varphi_\nu \circ \varphi_\mu^{-1})(x), H(\varphi_\nu \circ \varphi_\mu^{-1})(x) \big)
\end{equation}
the Jacobian and Hessian of the change of coordinates. The isomorphism between the bundle that we have just described and $\mathrm{Diff}_\mathrm{Sch}^{\, n} M$ is given by (notation as in \eqref{gluing}) $[\lambda, (t,x), (\sigma, \eta)] \mapsto [\bbF(x,t)]$, the class represented by any $\bbF(t,x)$ in the numerator of \eqref{schdiffbundle} s.t.\ $\bbF_\gamma^k = \sigma_\gamma^k$ for $\gamma = 1,\ldots, n$ and $\bbF_0^k + \sum_{\gamma = 1}^n \bbF_{\gamma\gamma}^k =  \eta^k$ w.r.t.\ the chart $\varphi_\lambda$.

There is a way of writing It\^o equations on a manifold so that all the coefficients, drift included, are vectors. It involves considering the additional structure of a linear connection $\nabla$ on $M$, i.e.\ a covariant derivative
\begin{equation}
	\nabla \colon TM \times \Gamma TM \to TM
\end{equation}
which is a smooth function that maps $T_xM \times \Gamma TM$ to $T_xM$, is $\bbR$-bilinear, and satisfies the Leibniz rule $\nabla_{U_x}(f V) = f(x) \nabla_{U_x}V + (U_xf)V_x$. Equivalently, a connection is described through its Hessian
\begin{equation}
	\nabla^2 \colon C^\infty M \to \Gamma (T^*M \otimes T^*M)
\end{equation}
which is an $\bbR$-linear map satisfying $\nabla^2(fg) = f \nabla^2 g + g \nabla^2 f + \dif f \otimes \dif g + \dif g \otimes \dif f$ for all $f, g \in C^\infty M$. These two data are equivalent and related by
\begin{equation}\label{covarianthessian}
	\langle \nabla^2_x f, V \otimes U \rangle = U_x(Vf) - (\nabla_{U_x}V) f
\end{equation}
If $\Gamma^{ij}_k$ are the Christoffel symbols of $\nabla$ w.r.t.\ a chart $\varphi$ (this means $\nabla_{\partial_x \varphi_i}\partial \varphi_j = \Gamma^k_{ij}(x) \partial_x \varphi_k$), the Hessian can be written as
\begin{equation}\label{hessiancoords}
	\nabla^2_x f = (\partial^2_x \varphi_{ij} - \Gamma_{ij}^k(x) \partial_x \varphi_k)(f) \dif_x \varphi^i \otimes \dif_x \varphi^j	
\end{equation} 
We will only be interested in connections modulo torsion, so it is not limiting for us to assume that a connection is symmetric or torsion-free, i.e.\ that its torsion tensor $\langle \tau_\nabla, U \otimes V \rangle = \nabla_U V - \nabla_V U - [U,V]$ vanishes, or equivalently that its Hessian is valued in $\Gamma(T^*M \odot T^*M)$. By far the most important example of such a connection is the Levi-Civita connection of a Riemannian metric $\mathscr g$; in this case the Hessian takes the form $\langle \nabla^2_x f, U_x \otimes V_x \rangle = \mathscr g( \nabla_{U_x} \mathrm{grad}^\mathscr{g}f, V_x )$. Torsion-free connections are relevant to our study of SDEs in that they correspond to the splittings of \eqref{SES}, i.e.\ a linear left inverse $\mathscr q$ to $\mathscr i$ or a linear right inverse $\mathscr j$ to $\mathscr p$  
\begin{equation}\label{splitSES}
	\begin{tikzcd}
		0 \arrow[r] &TM \arrow[r,"\mathscr i"] &\bbT M \arrow[r,swap,"\mathscr p"] \arrow[l,bend left, "\mathscr q"] &TM \odot TM \arrow[r] \arrow[l,bend right, swap, "\mathscr j"] & 0
	\end{tikzcd}
\end{equation}
The existence of the bundle maps $\mathscr j$ and $\mathscr q$ are equivalent to one another and to the the isomorphism $(\mathscr q, \mathscr p) \colon \bbT M \to TM \oplus (TM \odot TM)$ (this is the well-known splitting lemma \cite[p.147]{Hat02}, valid in the category of vector bundles). A torsion-free connection $\nabla$ on $M$ is equivalent to a splitting by setting
\begin{equation}\label{scriptq}
	(\mathscr q_x L_x)f \coloneqq L_xf - \langle \nabla^2_x f, \mathscr p_xL_x \rangle
\end{equation}
We recall that, given $V \in \Gamma TM$, $W_x \in T_xM$, the \say{composition} $U_x(V) \in \bbT_x M$ is defined by $U_x(V)f \coloneqq U_x(y \mapsto V_y f)$, and we have
\begin{equation}\label{pqofcomp}
	\mathscr p_x (U_x(V)) = U_x \odot V_x, \quad \mathscr q_x (U_x(V)) = \nabla_{U_x} V
\end{equation}
Using that $\partial^2_x \varphi_{ij} = \partial_x \varphi_i (\partial \varphi_j)$ and \eqref{hessiancoords} we have
\begin{equation}\label{qpartialij}
	\mathscr p_x \partial^2_x\varphi_{ij} = \partial_x \varphi_i \odot \partial_x \varphi_j, \quad \mathscr q_x \partial^2_x\varphi_{ij} = \Gamma_{ij}^k(x) \partial_x \varphi_k
\end{equation}
Another way to view this correspondence is by $\mathscr j^* \ddif_x f = \nabla^2_x f$.

Now, given symmetric connections on $N$ and $M$, a field of Schwartz morphisms $\bbF \in \Gamma\mathrm{Sch}(N,M)$ can be viewed as a field of block matrices
\begin{equation}\label{schwartzsplitting}
	\begin{bmatrix}
		F & G \\ 0 & F \otimes F
	\end{bmatrix}(x,z) \colon T_zN \oplus (T_zN \odot T_zN) \to T_xM \oplus (T_xM \odot T_xM)
\end{equation}
One can then require that $G \equiv 0$, so that $\bbF$ reduces to $F$, which defines the \emph{It\^o equation}
\begin{equation}\label{trueito}
	\dif X_t = F(X_t,Z_t) \dif Z_t
\end{equation}
Such equations have been considered in \cite{E90}. The data needed to define this equation is the same as that involved in the definition of the Stratonovich equation \eqref{StratonovichSDE}, namely an element of $\Gamma \mathrm{Hom}(TN,TM)$, but the meaning of the equation depends on the connections on $N$ and $M$. In local coordinates, using \eqref{qpartialij} to specify $\bbF_{\alpha\beta}^k$ in \eqref{SMincoords} to the case $G \equiv 0$, this equation takes the form
\begin{align}\label{itoeqcoords}
	\begin{split}
		\dif X_t^k &= F^k_\gamma(X_t,Z_t) \dif Z_t^\gamma \\
		&\mathrel{\hphantom{=}}+ \tfrac 12 \big( {^N\!}\Gamma_{\alpha \beta}^\gamma(Z_t) F^k_\gamma(X_t,Z_t)  - {^M\!}\Gamma_{ij}^k(X_t) F^i_\alpha F^j_\beta(X_t,Z_t) \big) \dif[Z^\alpha, Z^\beta]_t
	\end{split}
\end{align}
Note that if the Christoffel symbols on both manifolds vanish the above equation reduces to its first line; however, unless a manifold is flat a chart cannot in general be chosen so that the Christoffel symbols vanish (except for at a single chosen point: these are called normal coordinates). It\^o equations can be equivalently defined through the It\^o integral 
\begin{equation}\label{itoint}
	\int_0^\cdot \langle H_s, \dif X_s \rangle \coloneqq \int_0^\cdot \langle \mathscr q^*H_s, \ddif X_s \rangle
\end{equation}
by proceeding as in \eqref{fromintegraltoSDE}.

Recall that an $(M,\nabla)$-valued semimartingale is a \emph{local martingale} if for all $f \in C^\infty M$
\begin{equation}\label{martingaleproperty}
	f(X) - \int_0^\cdot \langle \mathscr p^* \nabla^2_{X_s} f , \ddif X_s \rangle
\end{equation}
is a real-valued local martingale (the integral is to be interpreted as half the quadratic variation of $X$ along the bilinear form $\nabla^2 f$); this property coincides with the usual local martingale property when $M$ is a vector space. In local coordinates an application of \eqref{hessiancoords} and \eqref{ddif} shows that the local martingale property corresponds to the requirement that
\begin{equation}\label{lmcoords}
	\dif X^k_t + \tfrac 12 \Gamma^k_{ij}(X_t) \dif[X^i,X^j]_t
\end{equation}
be a real-valued local martingale for each $k$. The It\^o integral \eqref{itoint} and It\^o equations \eqref{trueito} on manifolds behave well w.r.t.\ local martingales: if the integrand or driver is a local martingale, so is the integral or solution; this is again seen in local coordinates \eqref{itoeqcoords}.

In the following example we examine the case of diffusions, defined using It\^o equations, in which the issue of the drift not being a vector is (partially) resolved: 
\begin{expl}[It\^o diffusion]\label{pureitodiffusion}
	\autoref{SchMeydiffusion} specified to the above case ($M$ has a symmetric connection, $G \equiv 0$ in \eqref{schwartzsplitting}) becomes the equation
	\begin{equation}\label{pureitoeq}
		\dif X_t = \sigma_\gamma(X_t,t) \dif W^\gamma_t + \mu(X_t,t) \dif t
	\end{equation}
	where now $\mu(x,t) = \bbF(x,t) \in T_xM$ can legitimately be referred to as the \say{drift vector}. Note however that in an arbitrary chart $\varphi$ the drift will still carry a correction term:
	\begin{equation}\label{pureitocoordinates}
		\dif X_t^k = \sigma^k_\gamma(X_t,t) \dif W_t^\gamma + \bigg( \mu^k(X_t,t) -  \frac 12 \sum_{\gamma = 1}^n  \Gamma_{ij}^k(X_t) \sigma^i_\gamma \sigma^j_\gamma(X_t,t) \bigg) \dif t
	\end{equation}
	which reduces to the ordinary It\^o lemma if $M = \bbR^m$ and the chart $\varphi$ is a diffeomorphism of $\bbR^m$. The ${^N\!}\Gamma_{\alpha\beta}^\gamma$'s do not appear since the driver is already valued in a Euclidean space. The data needed to define such an equation coincides with that needed for \eqref{StratonovichSDE}, so we can define the bundle
	\begin{equation}
		\mathrm{Diff}_\mathrm{Ito}^{\, n} M \coloneqq \mathrm{Diff}_\mathrm{Strat}^{\, n} M \to M \times \bbR_{\geq 0}
	\end{equation}
	already defined in \eqref{stratdiffbundle}. Crucially, however, the Stratonovich and It\^o calculi give different meanings to the equation defined by a section of this bundle; in particular, a torsion-free connection on $M$ is required in the latter case. The \say{It\^o} and \say{Strat} therefore do not represent differences in the bundles, which are identical, but only serve as a reminder of which calculus is being used to give the section the meaning of an SDE.
\end{expl}
It\^o equations on manifolds are the true generalisation of their Euclidean space-valued counterparts, but have the disadvantage of only being defined w.r.t.\ a specific connection. For instance, if $F \in \Gamma \mathrm{Diff}^{\, n}_{\mathrm{Ito}}$, $M$ is Riemannian with $M'$ a Riemannian submanifold s.t. for all $z$ and $x \in M'$, $F(z,x)$ maps to $T_xM'$, $F$ does not in general define an It\^o equation on $M'$, since the Riemannian connection on $M'$ is not in general the restriction of that of $M$. However, $F$, seen as a field of Schwartz morphisms, does define a Schwartz-Meyer equation on $M'$ (with a $G$ term that is in general non-zero w.r.t.\ to the Riemannian connection on $M'$).

In the following table we summarise the advantages of these three ways of representing SDEs on manifolds:
\begin{center}
	\begin{tabular}{|c|c|c|c|}
		\hline
		& Stratonovich & Schwartz-Meyer/2-jet & It\^o  \\
		\hline
		Does not require $\nabla$ & \checkmark & \checkmark &  \\
		\hline
		Uses It\^o integration &  & \checkmark & \checkmark \\
		\hline
		Coefficients are vectors & \checkmark &  & \checkmark \\
		\hline
	\end{tabular}
\end{center}
It is natural to ask how these three types of equations are related to one another. In the case of diffusions, there exists a commutative diagram of bijections
\begin{equation}\label{ItoStratconvbdl}
	\begin{tikzcd}[column sep = tiny]
		& \Gamma \mathrm{Diff}_\mathrm{Sch}^{\, n} M \arrow[rd,rightarrow,"\mathscr b"] \\
		\Gamma \mathrm{Diff}_\mathrm{Strat}^{\, n} M \arrow[ur,rightarrow,"\mathscr a"] \arrow[rr,rightarrow,"\mathscr c"] && \Gamma \mathrm{Diff}_\mathrm{Ito}^{\, n} M
	\end{tikzcd}
\end{equation}
All three $\mathscr a,\mathscr b,\mathscr c$ are the identity on the diffusion coefficients. The behaviour of $\mathscr a, \mathscr b, \mathscr c$ on the Stratonovich, Schwartz-Meyer and It\^o drifts is explained below
\begin{equation}\label{abcformulae}
	\mathscr a  b  \coloneqq b + \frac 12 \sum_{\gamma = 1}^n \sigma_\gamma(\sigma_\gamma), \quad \mathscr b \eta  \coloneqq \mathscr q \eta, \quad \mathscr c  b  \coloneqq b + \frac 12 \sum_{\gamma = 1}^n \nabla_{\sigma_\gamma} \sigma_\gamma
\end{equation}
Note that, while $\mathscr b$ and $\mathscr c$ depend on the connection, $\mathscr a$ does not. If $\eta = \bbF_0 + \frac 12 \sum_{\gamma = 1}^n \bbF_{\gamma\gamma}$ is a Schwartz-Meyer drift, \eqref{Schcondition} and \eqref{pqofcomp} force $\eta - \frac 12\sum_{\gamma = 1}^n \sigma_\gamma(\sigma_\gamma)$ to lie in $T_x M$, which is thus $\mathscr a^{-1} \eta$. Moreover, we have $\mathscr b^{-1} \mu = \mathscr i \mu + \frac 12 \sum_{\gamma = 1}^n \mathscr j (\sigma_\gamma \odot \sigma_\gamma)$ and $\mathscr c^{-1} \mu = b - \frac 12 \sum_{\gamma = 1}^n \nabla_{\sigma_\gamma} \sigma_\gamma$. $\mathscr a, \mathscr b, \mathscr c$ define correspondences of SDEs in the sense that solutions are preserved (e.g.\ $X$ is a solution of $F \in \mathrm{Diff}_\mathrm{Strat}^{\, n} M$ if and only if $X$ is a solution of $\mathscr a F$, and the same for $\mathscr b, \mathscr c$). This is immediate by the expression of such equations in charts, by \eqref{pqofcomp} and the usual It\^o-Stratonovich conversion formula.
\begin{rem}
	What makes It\^o-Stratonovich conversion formulae difficult to state in the case of a general manifold-valued semimartingale driver $Z$, is that the change of calculus involves the emergence of new drivers which are not naturally valued in the manifold where Z is valued (the quadratic covariation of $Z$). Nevertheless, the map $\mathscr a$ can be defined in this general setting \cite[Lemma 7.22]{E89}, though its inverse cannot canonically.
\end{rem}

\section{Manifolds embedded in $\bbR^d$}\label{sec:projEmbedded}
In this chapter we will mostly be concerned with manifolds embedded in $\bbR^d$: these can be studied using the extrinsic, canonical, $\bbR^d$-coordinates instead of non-canonical local ones.
Let $M$ be an $m$-dimensional smooth manifold embedded in $\bbR^d$. We assume $M$ to be locally given by a non-degenerate Cartesian equation $F(x) = 0$: $M$ can be described globally in this way if and only if it is closed and its embedding has trivial normal bundle;
therefore, to preserve generality, we only assume $F$ to be local. Throughout this chapter the letter $x$ will denote a point in $\bbR^d$ and the letter $y$ a point in $M$. Thus $F \colon \bbR^d \to \bbR^{d-m}$ is a submersion, which implies $JF(x)JF(x)^\intercal \in GL(\bbR,d-m)$ for all $x \in \bbR^d$ ($JF(x) \in \bbR^{(d-m) \times d}$ the Jacobian of $F$ at $x$):
\begin{equation}
	JF(x)JF(x)^\intercal v^\intercal = 0 \ \Rightarrow \ (vJF(x))(vJF(x))^\intercal = v JF(x)JF(x)^\intercal v^\intercal = 0 \ \Rightarrow \ v = 0
\end{equation}
Let $\pi$, defined on a tubular neighbourhood $T$ of $M$ in $\bbR^d$ be the Riemannian submersion 
\begin{equation}\label{defpi}
	\pi(x) \coloneqq \arg\min \{ |x - y| : y \in M \}
\end{equation}
This map can be seen to exist by using the normal exponential map defined in \cite[p.132]{P06}, and is constant on the affine $(d-m)$-dimensional slices of $T$ which intersect $M$ orthogonally: this is because the fibre $\pi^{-1}(y)$ coincides with the union of all geodesics in $\bbR^d$ (i.e.\ straight line segments) which start at $y$, with initial velocity orthogonal to $M$, each taken for $t$ in some open interval containing $0$. 
It is important also to remember that $\pi$ is unique given the embedding of $M$ (on a thin enough $T$ such that it is well defined), whereas $F$ is not canonically determined. In what follows we will be concerned with understanding which quantities are dependent on the chosen $F$ and which instead only depend on the embedding of $M$. The only properties of $\pi$ that we will need are that
\begin{equation}\label{piprop}
	F \circ \pi \equiv 0, \quad \pi |_M = \mathbbm 1_M \Rightarrow \pi \circ \pi \equiv \pi
\end{equation}
Differentiating these (the second up to order 2) we obtain
\begin{align}\label{pisquared}
	\begin{split}
		\frac{\partial F}{\partial x^h}(\pi(x)) \frac{\partial \pi^h}{\partial x^k}(x) &= 0 \\
		\frac{\partial \pi}{\partial x^h}(\pi(x)) \frac{\partial \pi^h}{\partial x^k}(x) &= \frac{\partial \pi}{\partial x^k}(x) \\
		\frac{\partial^2 \pi}{\partial x^a \partial x^b}(\pi(x)) \frac{\partial \pi^a}{\partial x^i} \frac{\partial \pi^b}{\partial x^j}(x) + \frac{\partial \pi}{\partial x^h}(\pi(x)) \frac{\partial^2 \pi^h}{\partial x^i \partial x^j}(x)
		&= \frac{\partial^2 \pi}{\partial x^i \partial x^j}(x) 
	\end{split}
\end{align}
If $V_y \in T_M$ and $X$ is a smooth curve s.t.\ $X_0 = y$ and $\dot X_0 = V(y)$, differentiating $\pi(X_t) = X_t$ results in $J\pi(y) = V_y$: this shows that $J\pi|_{TM} = \mathbbm 1_{TM}$. By a similar argument, the fact that $\pi^{-1}(y)$ is a straight line segment that intersects $M$ orthogonally implies that $J\pi|_{T^\bot M} = \mathbbm 1_{T^\bot M}$ ($T^\bot_y M$ the normal bundle of $M$ at $y$). These two statements mean that
\begin{equation}\label{JPQ}
	P(y) = J\pi(y) \quad \text{for } y \in M
\end{equation}
where $P(y) \colon T_y \bbR^d \to T_yM$ is the orthogonal projection onto the tangent bundle of $M$, which can be defined in terms of $F$ as
\begin{align}\label{PQ}
	\begin{split}
		P(x) &\coloneqq \mathbbm 1 - Q(x) \qquad \text{where}\\
		Q(x) &\coloneqq JF^\intercal(x)(JF(x)JF^\intercal(x))^{-1}JF(x) \in \bbR^{d\times d} \qquad  \text{and we have}  \\
		PQ(x) &= 0 = QP(x), \quad QQ(x) = Q(x) = Q^\intercal(x), \quad PP(x) = P(x)= P^\intercal(x) 
	\end{split}
\end{align}
The notation is borrowed from \cite{CDL15}. Note that we can use $F$ to define $P,Q$ on a tubular neighbourhood of $M$, but these will only be independent of $F$ on $M$.
$Q(y)\colon T_y \bbR^d \to T_y^\bot M$ is the orthogonal projection onto the normal bundle. Another consequence of \eqref{pisquared} (evaluated at $y \in M$) that will be useful is that, for $V_y, W_y \in T_y\bbR^d$, and denoting $\overline U_y = P(y) U_y, \widecheck U_y = Q(y) U_y$
\begin{align}\label{linehat}
	\begin{split}
		&\mathrel{\hphantom{\Longrightarrow}} \frac{\partial^2 \pi}{\partial x^i \partial x^j}(y) \overline V^i_y \overline W^j_y \in T_y^\bot M, \quad \frac{\partial^2 \pi}{\partial x^i \partial x^j}(y) \overline V^i_y \widecheck W^j_y \in T_y M, \ \frac{\partial^2 \pi}{\partial x^i \partial x^j}(y)\widecheck V^i_y \widecheck W^j_y = 0 \\
		&\Longrightarrow  \frac{\partial^2 \pi}{\partial x^i \partial x^j}(y) V_y^i W_y^j = \underbrace{\frac{\partial^2 \pi}{\partial x^i \partial x^j}(y) \big(\overline V_y^i \overline W_y^j \big)}_{\in T^\bot_yM} + \underbrace{\frac{\partial^2 \pi}{\partial x^i \partial x^j}(y) \big(\overline V^i_y \widecheck W^j_y + \widecheck V^i_y \overline W^j_y \big)}_{\text{both terms } \in T_yM} 
	\end{split}
\end{align}
Actually, to show that the third term statement in the first line, we need a separate argument:
\begin{rem}\label{hessianonlydepends}
	Let $U \subseteq \bbR^d$, $f \colon U \to \bbR^e$, $y \in U$, $A_y,B_y \in T_y\bbR^d$. Then
	\begin{equation}
		\frac{\partial^2 f}{\partial x^i \partial x^j}(y) A^i_y B^j_y
	\end{equation}
	only depends on $f$ restricted to the affine plane (or line) centred in $y$ and spanned by $A_y,B_y$. Indeed, intending with $A$ the extension of $A_y$ to a constant vector field on $U$, we can write
	\begin{equation}
		\frac{\partial^2 f}{\partial x^i \partial x^j}(y) A^i_y B^j_y = \frac{\partial}{\partial x^j} \bigg|_y \bigg( \underbrace{\frac{\partial f}{\partial x^i}(x) A^i_x}_{\eqqcolon g(x)} \bigg) B^j_y
	\end{equation}
	This is the directional derivative of $g$ at $y$ in the direction $B_y$, and therefore only depends on the restriction of $g$ to the affine line $\mathrm{span}\{B_y\}$. But $g(x)$ is itself a directional derivative, and only depends on $f$ restricted to the affine line $\mathrm{span}\{ A_x \}$. Thus the whole expression only depends on $f$ restricted to $\bigcup_{x \in \mathrm{span}\{B_y\} } \mathrm{span}\{ A_{x}\} = \mathrm{span}\{ A_y, B_y\}$.
\end{rem}
This shows that the term in question only depends on $\pi$ restricted to $\mathrm{span}\{ \widecheck V_y, \widecheck W_y \}$, which is the constant $y$ map, whose derivatives therefore vanish.
\begin{rem}
	The other terms appearing in \eqref{linehat} have a description that should be more familiar to differential geometers: 
	\begin{align}\label{nablad2pi}
		\begin{split}
			\frac{\partial^2 \pi}{\partial x^i \partial x^j} (y) \overline V^i_y \overline W^j_y &= {^{\bbR^d}\nabla_{\overline V_y}^\bot} \overline W \coloneqq Q(y){^{\bbR^d}\nabla_{\overline V_y}} \overline W =  \mathrm{I\!I} \big(\overline V_y,\overline W_y \big)	\\
			-\frac{\partial^2 \pi}{\partial x^i \partial x^j} (y)  \overline V^i_y \widecheck W^j_y &= {^{\bbR^d}\nabla_{\overline V_y}^\top} \widecheck W \coloneqq P(y){{^{\bbR^d}}\nabla_{\overline V_y}} \widecheck W
		\end{split}
	\end{align}
	where ${^{\bbR^d}}\nabla$ denotes covariant differentiation in $\bbR^d$ (i.e.\ just directional differentiation). Notice this is true independently of the chosen extension of $\overline W, \widecheck W$ to local vector fields, a priori needed to give the RHSs a meaning. The first term is the second fundamental form of $\overline V_y, \overline W_y$ \cite[p.134]{L97}, whereas the second term is the second fundamental tensor \cite[Def.\ 3.6.1]{J05}. If $M$ is an open set of an affine subspace of $M$, $\pi$ is a linear map and both terms vanish. We prove the first of the two equalities in \eqref{nablad2pi}, the second is proved similarly:
	\begin{equation}\label{similarargument}
		Q(y){^{\bbR^d}\nabla_{\overline V_y}} \overline W = Q_j(y) \frac{\partial \overline W^j}{\partial x^i}(y) \overline V^i_y = -\frac{\partial Q_j}{\partial x^i}(y) \overline W^j_y \overline V^i_y = \frac{\partial^2 \pi}{\partial x^i \partial x^j}(y)\overline V^i_y \overline W^j_y
	\end{equation}
	where the second equality follows from the fact that $Q \overline W = 0$ (and that the derivative is taken in a tangential direction, i.e.\ $\overline V_y \in T_y M$), and the last equality is given by \eqref{dP} below. Note that the terms of \eqref{nablad2pi} are extrinsic, in the sense that they depend on the embedding of $M$, unlike
	\begin{equation}\label{covariantembedded}
		{^M \!}\nabla_{\overline V_y} \overline W_y = P(y) {^{\bbR^d}\nabla_{\overline V_y}} \overline W
	\end{equation}
	the Levi-Civita connection of the Riemannian metric on $M$, which is intrinsic to $M$.
\end{rem}

Finally, it will be necessary to consider the relationship between the derivatives of $P,Q$ and the second derivatives of $\pi$. We differentiate \eqref{JPQ} at time $0$ along a smooth curve $Y_t$ in $M$ with $Y_0 = 0$ and $\dot Y_0 = \overline V_y \in T_y M$ and obtain
\begin{equation}\label{eq:somethingPpi}
	\frac{\partial P_k}{\partial x^h}(y) \overline V^h_y = \frac{\partial^2 \pi}{\partial x^i \partial x^j}(y)  \overline V^i_y
\end{equation}
from which we obtain, for $W \in T_yM$
\begin{align}\label{dP}
	\begin{split}
		-\frac{\partial Q_k}{\partial x^h}(y) \overline V^h_y \overline W^k_y = \frac{\partial P_k}{\partial x^h}(y) \overline V^h_y \overline W^k_y &= \frac{\partial^2 \pi^k}{\partial x^i \partial x^j}(y) \overline V^i_y \overline W_y^j \in T_y^\bot M \\
		-\frac{\partial Q_k}{\partial x^h}(y) \overline V^h_y \widecheck W^k_y = \frac{\partial P_k}{\partial x^h}(y) \overline V^h_y \widecheck W^k_y &= \frac{\partial^2 \pi}{\partial x^i \partial x^j}(y)  \overline V^i_y \widecheck W^j_y \in T_yM
	\end{split}
\end{align}
where we have used \eqref{linehat}.\\

We now consider a setup $\mathscr S = (\Omega, \mathcal F, P)$ satisfying the usual conditions, $W$ an $n$-dimensional Brownian motion defined on $\mathscr S$. Consider the $W$-driven diffusion Stratonovich SDE
\begin{equation}\label{StratSDE}
	\dif X_t^k = \sigma_\gamma^k(X_t,t) \circ \dif W_t^\gamma + b^k(X_t,t) \dif t, \quad X_0 = y_0 \in M
\end{equation}
As already discussed in \autoref{sec:SDEs}, the natural condition on $\sigma_\gamma, b$ which guarantees that $X$ will stay on $M$ for its lifetime is their tangency to $M$:
\begin{equation}\label{StratonovichTang}
	Q(y) \sigma_\gamma (y,t) = 0 = Q(y) b(y,t) \quad \text{for all } y \in M, t \geq 0, \gamma = 1,\ldots, n
\end{equation} 
Our focus, however, will be mostly on the It\^o SDE
\begin{equation} \label{sde}
	\dif X_t^k = \sigma^k_\gamma(X_t,t) \dif W^\gamma_t + \mu^k(X_t,t) \dif t, \quad X_0 = y_0 \in M
\end{equation}
with smooth coefficients defined in $[0,+\infty) \times \bbR^d$; we do not assume them to be globally Lipschitz, so the solution might only exist up to a positive stopping time, not in general bounded from below by a positive deterministic constant. We are interested in deriving the \say{tangency condition} for the above SDE, i.e.\ a condition on the coefficients that will guarantee that the solution will not leave $M$. One way to impose this is to convert \eqref{sde} to Stratonovich form
\begin{equation}\label{itotostratonovich}
	\dif X_t^k = \sigma^k_\gamma(X_t,t) \circ \dif W^\gamma_t + \bigg(\mu^k - \frac 12 \sum_{\gamma = 1}^n \sigma^h_\gamma \frac{\partial \sigma^k_\gamma}{\partial x^h} \bigg)(X_t,t) \dif t, \quad X_0 = y_0 \in M
\end{equation}
and require \eqref{StratonovichTang}: 
\begin{equation} \label{strattan}
	\begin{dcases}
		Q_k(y) \sigma^k_\gamma(y,t) = 0 \\ Q_k(y) \bigg( \mu^k - \frac 12 \sum_{\gamma = 1}^n \sigma^h_\gamma \frac{\partial \sigma^k_\gamma}{\partial x^h}   \bigg)(y,t) = 0
	\end{dcases}
\end{equation}
Now, given that $Q \sigma_\alpha$ vanishes on $M$, all its directional derivatives along the tangent directions $\sigma_\beta$ will too, which gives, using \eqref{dP}
\begin{equation}\label{itotrick}
	0 = \frac{\partial (Q \sigma_\alpha)}{\partial x^h} \sigma_\beta^h = \frac{\partial Q_i}{\partial x^j} \sigma^i_\alpha \sigma^j_\beta + Q_k \frac{\partial \sigma^k_\alpha}{\partial x^h}  \sigma^h_\beta \Longrightarrow  Q_k  \frac{\partial \sigma^k_\alpha}{\partial x^h} \sigma^h_\beta = \frac{\partial^2 \pi}{\partial x^i \partial x^j} \sigma_\alpha^i \sigma_\beta^j \quad \text{on $M$}
\end{equation}
We can thus reformulate the second equation in \eqref{strattan} to obtain
\begin{equation} \label{itotan}
	\begin{dcases}
		Q_k(y) \sigma^k_\gamma(y,t) = 0 \\ Q_k(y) \mu^k(y,t) = \frac 12 \sum_{\gamma = 1}^n \frac{\partial^2 \pi}{\partial x^i \partial x^j}(y) \sigma^i_\gamma \sigma^j_\gamma(y,t)
	\end{dcases}
\end{equation}
This is useful because it removes the reliance of this constraint on the derivatives of $\sigma$, and can be interpreted as saying that the diffusion coefficients must be tangent to $M$ and the It\^o drift must instead lie on the space parallel to the tangent space of $M$, displaced by an amount which depends on the second fundamental form of $M$ applied to the diffusion coefficients.
\begin{rem}[Tangency of a second-order differential operator]\label{remarktangency}
	\eqref{itotan} can also be derived by writing the second order tangency condition for $L_y^k \partial_y x_k + L^{ij}_y \partial^2_y x_{ij} = L_y \in \bbT_y \bbR^d$ to belong to $\bbT_y M$: this is done by writing $\bbT_y \pi L_y = L_y$ in $\bbR^d$-coordinates as
	\begin{equation}\label{diffusortangency}
		\begin{bmatrix}
			L^h_y \\ L^{ab}_y
		\end{bmatrix} =
		\begin{bmatrix}
			\frac{\partial \pi^h}{\partial x^k} & \frac{\partial^2 \pi^h}{\partial x^i \partial x^j} \\
			0 & \frac{\partial \pi^a}{\partial x^i} \frac{\partial \pi^b}{\partial x^j}
		\end{bmatrix}(y) \begin{bmatrix}
			L^k_y \\ L^{ij}_y
		\end{bmatrix}
	\end{equation} 
	and then applying it to $L_y = \sigma_\gamma(y,t), \eta(y,t)$, given in terms a field of Schwartz morphisms $\bbF$ as
	\begin{equation}\label{SMembeddedSDE}
		\sigma_\gamma^k = \bbF_\gamma^k, \quad      \eta^k = \bbF_0^k + \frac 12 \sum_{\gamma = 1}^n \bbF^k_{\gamma\gamma}
	\end{equation}
	Note that it would instead be incorrect to split $\bbF$ according to the Euclidean connection into a matrix with $F$ and $G$ terms as in \eqref{schwartzsplitting}, and then to require that $F$ and $G$ map to $TM$, since the splitting of $\bbF$ according to the connection on $M$ will be different, i.e.\ the diagram
	\begin{equation}
		\begin{tikzcd}
			\bbT \bbR^d \arrow[r,leftrightarrow,"\cong"] & T \bbR^d \oplus (T\bbR^d \odot T\bbR^d)\\
			\bbT M \arrow[r,leftrightarrow,"\cong"] \arrow[u,hookrightarrow] 
			& TM \oplus (TM \odot TM) \arrow[u,hookrightarrow]
		\end{tikzcd}
	\end{equation} 
	does not commute.
\end{rem}

We now compute the Hessian for embedded $M$: for $f \in C^\infty M$ we have
\begin{equation}\label{hessian embedded}
	\big\langle {^M \!} \nabla^2_y f, \overline V_y \otimes \overline U_y \big\rangle = \big\langle {^{\bbR^d} \!} \nabla^2_y (f \circ \pi), \overline V_y \otimes \overline U_y \big\rangle 
\end{equation}
where we have used \eqref{covarianthessian}, \eqref{covariantembedded} to reduce this to a computation of directional derivatives, and finally \eqref{dP} (the argument is similar to \eqref{similarargument}). ${^{\bbR^d} \!} \nabla^2$ of course is just the ordinary Hessian. We can now compute ${^M \!\!} \mathscr q$, the splitting appearing in \eqref{splitSES} w.r.t.\ the connection on $M$: if $\bbT_yM \ni L_y = L_y^k \partial_y x_k + L^{ij}_y \partial^2_y x_{ij}$, using \eqref{scriptq} yields
\begin{align}
	\begin{split}
		({^M \!\!} \mathscr q_y L_y) f &=  L_y(f) - \big\langle {^M \!}\nabla^2_y f, \mathscr p_y L_y  \big\rangle \\
		&= L_y(f \circ \pi) - \big\langle {^{\bbR^d} \!}\nabla^2_y (f \circ \pi), L_y^{ij} \partial^2_y x_{ij} \big\rangle \\
		&= \frac{\partial f}{\partial x^h}(y) \frac{\partial \pi^h}{\partial x^k}(y) L^k_y 
	\end{split}
\end{align}
which means
\begin{equation}\label{qambient}
	{^M \!\!}\mathscr q_y = P(y) \circ {^{\bbR^d} \! \!} \mathscr q_y \colon \bbT_y M \to T_yM
\end{equation} 
Therefore the condition on an arbitrary Schwartz morphism of being It\^o w.r.t.\ to the Riemannian connection on $M$ in the sense of \autoref{pureitodiffusion} is ${^M \!\!}\mathscr q \circ \bbF \circ {^{\bbR^d} \!\!}\mathscr j = 0$, or ${^M \!\!}\mathscr q \bbF_{\alpha\beta} = 0$, which in $\bbR^d$-coordinates is
\begin{equation}\label{itoMcondition}
	P_k(y) \bbF_{\alpha\beta}^k(y,t) = 0
\end{equation}
Compare this with the stronger condition of $\bbF$ of being It\^o w.r.t.\ to the connection on $\bbR^d$, which is $\bbF_{\alpha\beta}^k(y,t) = 0$. Thus, given an It\^o equation $\bbF$ on $M$, defined as in \eqref{pureitoeq} ($\sigma_\gamma = \bbF_\gamma$, $\overline\mu = \bbF_0$) we have that the drift in $\bbR^d$ of such equation is given by $\overline\mu^k + \frac 12 \sum_{\gamma = 1}^n \bbF^k_{\gamma\gamma}$, with the first term tangent to $M$ and the second orthogonal to $M$, and equal to $\frac 12 \sum_{\gamma = 1}^n \frac{\partial \pi^h}{\partial x^i \partial x^j} \sigma^i \sigma^j_\gamma$, by \autoref{remarktangency} and \eqref{itoMcondition}. Therefore an It\^o equation on $M$ with coefficients $\sigma_\gamma, \overline\mu$ is read in ambient coordinates as
\begin{equation}\label{itocoefficientsembedded}
	\dif X_t^k = \sigma_\gamma^k(Y_t,t) \dif W_t^\gamma + \bigg(\overline\mu^k + \frac 12 \sum_{\gamma = 1}^n \frac{\partial^2 \pi^k}{\partial x^i \partial x^j} \sigma^i_\gamma \sigma^j_\gamma \bigg)(Y_t,t) \dif t
\end{equation}
Notice that the tangential part of the $\bbR^d$-drift, $\overline\mu$, is arbitrary, while its orthogonal part is determined by the diffusion coefficients, and the condition that the solution remain on $M$.
\begin{figure}[h]
	\minipage{0.5\textwidth}
	\includegraphics[width=\linewidth]{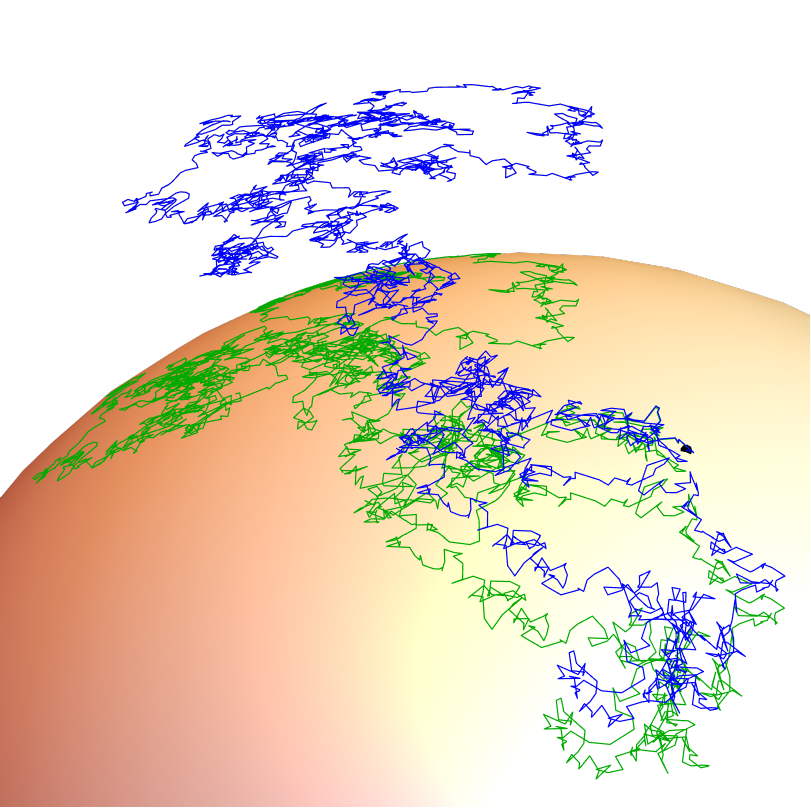}
	\vspace{-20pt}
	\endminipage\hfill
	\minipage{0.5\textwidth}
	\includegraphics[width=\linewidth]{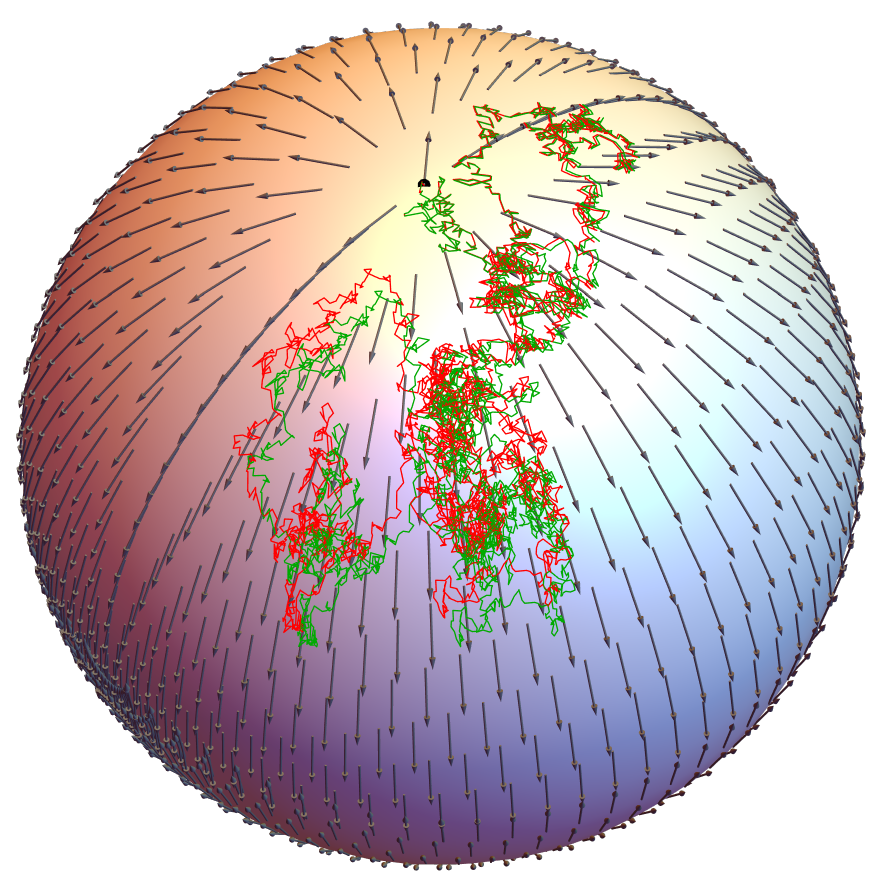}
	\vspace{-20pt}
	\endminipage
	\caption{On the left a sample path of the solution to the It\^o equation (blue) with the two diffusion coefficients $2(x^2+y^2+z^2)^{-1}(-y,x,0)$, $2(x^2+y^2+z^2)^{-1}(0,-z,y)$, which are tangent to $S^2 \hookrightarrow \bbR^3$, zero drift and initial condition $(0,1,0)$; in the same plot a sample path (using the same random seed) of the solution to the Stratonovich equation (green) defined by the same vector fields and initial condition. The solution to the It\^o equation drifts radially outwards, while the solution to the Stratonovich equation remains on $S^2$. On the right we compare the same Stratonovich path with a sample path of the solution to the It\^o equation (red) with the same diffusion coefficients and initial condition, but with the orthogonal drift term necessary to keep the solution on $S^2$ \eqref{itotan}. The resulting solution is an $S^2$-valued local martingale, while the solution to the Stratonovich equation is not: this is illustrated by plotting the vector field on $S^2$ given by tangential component of the It\^o drift possessed by the Stratonovich equation: this can be viewed as a manifold-valued drift component.}\label{fig:itoStrat}
\end{figure}
The notion of $M$-valued local martingale also has a description in terms of ambient coordinates \cite[\P 4.10]{E89}: for an $M$-valued It\^o process (such as the solution to \eqref{itocoefficientsembedded}) the local martingale property is equivalent to requiring that the drift be orthogonal to $M$ at each point (and thus determined by the diffusion coefficients; for \eqref{itocoefficientsembedded} this means $\overline \mu = 0$). This condition is very reminiscent of the property of geodesics of having acceleration orthogonal to $M$ \cite[Lemma 8.5]{L97}.

Using all \eqref{abcformulae} and \eqref{qambient} it is easy to verify that converting between Stratonovich, Schwartz-Meyer and It\^o equations on $M$ is equivalent when treating the equations as being valued in $M$ or in $\bbR^d$. By this we mean that, denoting with $\mathrm{Diff}_{\mathrm{Strat},M}^{\, n} \bbR^d$ the bundle of Stratonovich equations on $\bbR^d$ which restrict to equations on $M$ (and analogously for the other two diffusion bundles) the maps $\mathscr a, \mathscr b, \mathscr c$ of \eqref{ItoStratconvbdl} fit into the commutative diagram  
\begin{equation}\label{commutativePrism}
	\begin{tikzcd}[column sep = tiny]
		& \Gamma \mathrm{Diff}_{\mathrm{Sch},M}^{\, n} \bbR^d \arrow[rd,leftrightarrow,"{^{\bbR^{\scaleto{d}{3pt}}} \! \!}\mathscr b"] \arrow[from = ddd,leftarrow]  \\
		\Gamma \mathrm{Diff}_{\mathrm{Strat},M}^{\, n} \bbR^d \arrow[ur,leftrightarrow,"{^{\bbR^{\scaleto{d}{3pt}}} \! \! \! }\mathscr a"] \arrow[rr,leftrightarrow,"{^{\bbR^{\scaleto{d}{3pt}}} \! \! \!}\mathscr c",near end,crossing over] && \Gamma \mathrm{Diff}_{\mathrm{Ito},M}^{\, n} \bbR^d \\ \\
		& \Gamma \mathrm{Diff}_\mathrm{Sch}^{\, n} M \arrow[rd,leftrightarrow,"{^M \!}\mathscr b"] \\
		\Gamma \mathrm{Diff}_\mathrm{Strat}^{\, n} M  \arrow[uuu,leftarrow]\arrow[ur,leftrightarrow,"{^M \! \!}\mathscr a"] \arrow[rr,leftrightarrow,"{^M \! \!}\mathscr c"] && \Gamma \mathrm{Diff}_\mathrm{Ito}^{\, n} M \arrow[uuu,leftarrow]
	\end{tikzcd}
\end{equation}
where vertical arrows denote restriction. An embedding argument immediately allows us to extend this assertion to the case where $\bbR^d$ is substituted with a Riemannian manifold of which $M$ is a Riemannian submanifold. This confirms there is no ambiguity in converting an $M$-valued SDE between its various forms.
\begin{expl}[Time dependent submanifold]\label{timedepM}
	Observe that the tangency conditions \eqref{StratonovichTang} and \eqref{itotan} can be written respectively as
	\begin{equation}
		\begin{cases}
			(\mathbbm 1 - J\widetilde \pi) \sigma_\gamma = 0 \\
			(\mathbbm 1 - J\widetilde\pi)b = 0
		\end{cases}
		\quad \begin{dcases}
			(\mathbbm 1 - J\widetilde\pi) \sigma_\gamma = 0 \\
			(\mathbbm 1 - J\widetilde\pi)\mu = \frac 12 \sum_{\gamma = 1}^n \frac{\partial^2 \widetilde\pi}{\partial x^i \partial x^j}(y) \sigma^i_\gamma \sigma^j_\gamma
		\end{dcases}
	\end{equation}
	for \emph{any} smooth map $\widetilde\pi$ defined on a tubular neighbourhood of $M$, with values in $M$, s.t.\ $\widetilde\pi|_M = \mathbbm 1_M$, by the same exact reasoning (for the It\^o case we argue as in \autoref{remarktangency}). $J \widetilde \pi(y)$ is no longer the orthogonal projection $P(y)$, but still restricts to the identity on $T_yM$ for $y \in M$, i.e.\ it has the property that $\ker(\mathbbm 1 - J \widetilde \pi) = TM$ on $M$. Allowing ourselves to consider all such tubular neighbourhood projections is useful in the following application. Given that we are considering time-dependent equations, it is very natural to also allow the submanifold $M$ to be time-dependent. Making this precise entails considering a smooth $(m+1)$-dimensional manifold $\widetilde M$ embedded in $\bbR^{1+d}$, s.t.\ $M_t \coloneqq \widetilde M \cap \{x_0 = t\}$ is a smooth $m$-dimensional manifold embedded in $\{x_0 = t\} \times \bbR^d$. We are looking for conditions on $\sigma, b$ (resp. $\mu$) which are sufficient to guarantee the solution to \eqref{StratSDE} (resp. \eqref{sde}) $X_t$ to belong to $M_t$ for all $t$ for which it is defined. We then consider the $\bbR^{1+d}$-valued process $(t,X_t)$, which satisfies the dynamics
	\begin{equation}\label{tX}
		\dif \begin{bmatrix}
			t \\ X_t
		\end{bmatrix} = \begin{bmatrix}
			0 \\ \sigma(X_t,t)
		\end{bmatrix} \circ \dif W_t + \begin{bmatrix}
			1 \\ b(X_t,t)
		\end{bmatrix} \dif t \quad \text{resp.} \quad =\begin{bmatrix}
			0 \\ \sigma(X_t,t)
		\end{bmatrix} \dif W_t + \begin{bmatrix}
			1 \\ \mu(X_t,t)
		\end{bmatrix} \dif t
	\end{equation}
	Then, given a thin enough tubular neighbourhood of $\widetilde M$ in $\bbR^{1+d}$ consider the map 
	\begin{equation}
		\widetilde \pi \colon \widetilde T \to \widetilde M,\quad  \widetilde \pi (t,x) = \pi_t(x)
	\end{equation}
	where $\pi_t$ is defined as in \eqref{defpi} for the manifold $M_t$. Notice that this does not in general coincide with the Riemannian projection of a tubular neighbourhood onto $\widetilde M$, which in general has no reason to preserve time, i.e.\ be expressible as a union of $\pi_t$'s. The identity $J \widetilde \pi J \widetilde \pi = J \widetilde \pi$ can be written in block matrix form as
	\begin{equation}
		\left[
		\begin{array}{c|c c c}
			1 & 0 & \cdots & 0 \\
			\hline
			\\
			J\pi_t \dot\pi_t + \dot \pi_t &  & J\pi_t J\pi_t
			& 
			\\
			\phantom{\pi}
		\end{array}
		\right] = \left[
		\begin{array}{c|c c c}
			1 & 0 & \cdots & 0 \\
			\hline
			\\
			\dot\pi_t &  & J\pi_t
			& 
			\\
			\phantom{\pi}
		\end{array}
		\right]
	\end{equation}
	where we are denoting $\dot \pi_t(y) = \frac{\dif}{\dif t} \pi_t(y)$: this implies that at each point $y \in M_t$, $\dot\pi_t(y) \in T_y^\bot M_t$. This choice of the tubular neighbourhood projection will be further motivated later on, in \autoref{tproj}, \autoref{optMt}. In view of the above considerations, we can use it anyway to impose tangency of the SDE: this results in an unmodified condition on the diffusion coefficients, and the conditions on the orthogonal components of the Stratonovich and It\^o drifts are given respectively by
	\begin{align}
		\begin{split}
			(\mathbbm 1 - J \pi_t) b(y,t) &= \dot \pi_t(y) \\
			(\mathbbm 1 - J \pi_t) \mu(y,t) &= \frac 12 \sum_{\gamma = 1}^n \frac{\partial^2 \pi_t}{\partial x^i \partial x^j}(y) \sigma^i_\gamma \sigma^j_\gamma + \dot \pi_t(y)
		\end{split}
	\end{align}
	which keep track of the evolution of $M_t$ in time.
\end{expl}

\section{Projecting SDEs}\label{sec:projecting}
In \autoref{sec:SDEs} we discussed three ways of representing SDEs on manifolds: Stratonovich, Schwartz-Meyer and It\^o. In this section we will define, for each one of these representations, a natural projection of the SDE onto a submanifold. We will mostly take the ambient manifold to be $\bbR^d$, which will allow us to use the theory of the previous section to derive formulae for the projections in ambient coordinates.

Let $M$ be a smooth submanifold of the smooth manifold $D$, let $T$ be a tubular neighbourhood of $M$ in $D$ and 
\begin{equation}\label{pigeneral}
	\text{$\pi \colon T \to M$ a smooth map which restricts to the identity on $M$}
\end{equation}
If $D$ is Riemannian $\pi$ can be chosen as in \eqref{defpi}, but this is not necessary. Let $F \in \Gamma \mathrm{Hom}(TN, TD)$ be a Stratonovich equation driven by an $N$-valued semimartingale $Z$, where $N$ is another smooth manifold. We can then define the $M$-valued Stratonovich equation 
\begin{equation}\label{GeneralStratProj}
	M \times N \ni (y,z) \mapsto \widetilde F(y,z) \coloneqq T_y\pi \circ F(y,z) \in \mathrm{Hom}(T_zN, T_yM)
\end{equation}
We call this Stratonovich SDE the \emph{Stratonovich projection of $F$}.

Now consider the $Z$-driven, $D$-valued Schwartz-Meyer equation $\bbF \in \Gamma\mathrm{Sch}(N,M)$. We can project this SDE to an SDE on $M$ too, by
\begin{equation}\label{GeneralSMProj}
	M \times N \ni (y,z) \mapsto \widehat \bbF(y,z) \coloneqq \bbT_y\pi \circ \bbF(y,z) \in \mathrm{Sch}_{z,y}(N,M)
\end{equation}
We call this Schwartz-Meyer SDE the \emph{It\^o-jet projection of $\bbF$}.

If $N$, $D$ and $M$ all carry torsion-free connections we can interpret a section $F \in \Gamma\mathrm{Hom}(TN, TD)$ as an It\^o equation, and similarly for
\begin{equation}\label{GeneralItoProj}
	M \times N \ni (y,z)  \mapsto \overrightarrow F(y,z) \coloneqq T_y\pi \circ F(y,z) \in \mathrm{Hom}(T_zN, T_yM)
\end{equation}
We call this It\^o SDE the \emph{It\^o-vector projection} of $F$. Most often $D$ will be Riemannian, so that Levi-Civita connections are defined on both $D$ and $M$. Note that the It\^o-vector projection is identical to the Stratonovich projection as a map, but the interpretations of the resulting sections as SDEs differ (and the It\^o-vector projected SDE depends explicitly on the connections on all three manifolds). The names of these three projections are taken from \cite{ABR19}, where they were first defined.
\begin{rem}[Naturality of the SDE projections]
	Assume we have a commutative square
	\begin{equation}
		\begin{tikzcd}
			D \arrow[r,"\phi"] \arrow[d,"\pi"] & D' \arrow[d,"\pi'"] \\
			M \arrow[r,"\phi|_M"] & M'
		\end{tikzcd}
	\end{equation}
	where $\phi$ a diffeomorphism, $D,M,\pi$ as above, and similarly for $D',M',\pi'$. Then functoriality of $T$ and $\bbT$ imply that the Stratonovich and It\^o-jet projections are natural in the sense that the squares
	\begin{equation}
		\begin{tikzcd}
			\mathrm{Hom}(TN,TD) \arrow[r,"T\phi"] \arrow[d,"\widetilde{\phantom{aa}}"] & \mathrm{Hom}(TN,TD') \arrow[d,"\widetilde{\phantom{aa}}"] & \mathrm{Sch}(N,D) \arrow[r,"\bbT \phi"] \arrow[d,"\widehat{\phantom{aa}}"] & \mathrm{Sch}(N,D') \arrow[d,"\widehat{\phantom{aa}}"] \\
			\mathrm{Hom}(TN,TM) \arrow[r,"T\phi|_M"] & \mathrm{Hom}(TN,TM') & \mathrm{Sch}(N,M) \arrow[r,"\bbT \phi|_M"] & \mathrm{Sch}(N,M')
		\end{tikzcd}
	\end{equation}
	commute. The It\^o-vector projection cannot be natural in the same way, since we are still free to modify the connections on all four manifolds. However, if $D,D'$ are Riemannian and $\phi$ is a global isometry, the corresponding statement does hold for the It\^o-vector projection as well: this is by naturality of the Levi-Civita connection \cite[Proposition 5.6]{L97}.
\end{rem}
\begin{rem}[The It\^o-vector projection preserves local martingales]
	Although the It\^o-vector projection is natural w.r.t.\ a smaller class of maps, it has the advantage of preserving the local martingale property: by this we mean that if the driver is a local martingale, so must the solution to the It\^o-vector-projected SDE be. This is shown simply by the good behaviour of It\^o equations w.r.t.\ manifold-valued local martingales.
\end{rem}
\begin{rem}
	One might wonder whether it is possible to \say{push forward} SDEs according to an arbitrary smooth and surjective map $f \colon D \to D'$. If $f$ is a surjective function admitting a smooth right inverse $\iota$, then we may write the pushforward of, say, the Stratonovich SDE $\dif X = F(X,Z) \circ \dif Z$ as $\dif Y = F(Z,\iota(Y)) \circ \dif Y$. This condition on $f$ essentially corresponds to the condition \eqref{pigeneral}. For general smooth surjective $f$ (such as the bundle projection of a non-trivial principal bundle), however, we do not see a way of defining a new closed form SDE on $D'$.
\end{rem}

We will now restrict our attention to the projections of $\bbR^d$-valued diffusions onto the embedded manifold $M$. Focusing on diffusions has the advantage of allowing us to use the maps \eqref{ItoStratconvbdl} to compare the projections. In other words we can ask if the vertical rectangles in the diagram
\begin{equation}\label{NonCommutativePrism}
	\begin{tikzcd}[column sep = tiny]
		& \Gamma \mathrm{Diff}_\mathrm{Sch}^{\, n} \bbR^d \arrow[rd,leftrightarrow,"{^{\bbR^{\scaleto{d}{3pt}}} \! \!}\mathscr b"] \arrow[ddd,twoheadrightarrow,"\widehat{\phantom{aa}}", near end]  \\
		\Gamma \mathrm{Diff}_\mathrm{Strat}^{\, n} \bbR^d \arrow[ur,leftrightarrow,"{^{\bbR^{\scaleto{d}{3pt}}} \! \! \! }\mathscr a"] \arrow[ddd,twoheadrightarrow, "\widetilde{\phantom{aa}}"]\arrow[rr,leftrightarrow,"{^{\bbR^{\scaleto{d}{3pt}}} \! \! \!}\mathscr c",near end,crossing over] && \Gamma \mathrm{Diff}_\mathrm{Ito}^{\, n} \bbR^d \arrow[ddd,twoheadrightarrow,"\overrightarrow{\phantom{aa}}"] \\ \\
		& \Gamma \mathrm{Diff}_\mathrm{Sch}^{\, n} M \arrow[rd,leftrightarrow,"{^M \!}\mathscr b"] \\
		\Gamma \mathrm{Diff}_\mathrm{Strat}^{\, n} M  \arrow[ur,leftrightarrow,"{^M \! \!}\mathscr a"] \arrow[rr,leftrightarrow,"{^M \! \!}\mathscr c"] && \Gamma \mathrm{Diff}_\mathrm{Ito}^{\, n} M 
	\end{tikzcd}
\end{equation}
commute (compare with \eqref{commutativePrism}, in which the equations on top already restrict to equations on $M$). We will show that they do not, and that all combinations of possibilities regarding their non-commutativity are possible. Examples of these cases are to be found in \autoref{subsec:comparison} below. We recall the notation $\overline V_y \coloneqq P(y)V_y$, $\widecheck V_y \coloneqq Q(y)V_y$ and begin by considering the $\bbR^d$-valued Stratonovich SDE \eqref{StratSDE}. By \eqref{JPQ} the coefficients of the Stratonovich projection of this SDE will just be the projected coefficients: $\widetilde \sigma_\gamma = \overline \sigma_\gamma, \widetilde b = \overline b$, so that the resulting Stratonovich equation is
\begin{align}\label{ProjStratSDE}
	\begin{split}
		\dif Y_t &= \overline \sigma_\gamma(Y_t,t) \circ \dif W_t^\gamma + \overline b(Y_t,t) \dif t, \quad Y_0 = y_0 \in M \\
		&= \frac{\partial \pi}{\partial x^k}(Y_t) \sigma_\gamma^k(Y_t,t) \circ \dif W^\gamma_t + \frac{\partial \pi}{\partial x^k}(Y_t) b^k(Y_t,t) \dif t
	\end{split}
\end{align}
Throughout this chapter we will use $X$ for the initial SDE and $Y$ to denote the projected SDE. Now assume we start with \eqref{sde}, and want an It\^o SDE on $M$. We can still use the Stratonovich projection by converting the SDE to Stratonovich form as in \eqref{itotostratonovich}, projecting as above, and converting back to It\^o form (by \eqref{commutativePrism} this last conversion can be seen to occur interchangeably in $M$ or in $\bbR^d$). We have
\begin{align}
	\begin{split}
		\dif Y_t &= \overline \sigma_\gamma(Y_t,t) \circ \dif W^\gamma_t + P_k(Y_t) \bigg( \mu^k - \frac 12 \sum_{\gamma = 1}^n \sigma^h_\gamma \frac{\partial \sigma_\gamma^k}{\partial x^h}  \bigg)(Y_t,t) \dif t \\
		&= \overline \sigma_\gamma(Y_t,t) \dif W^\gamma_t + \bigg(\underbrace{\overline\mu + \frac 12 \sum_{\gamma = 1}^n \bigg(\overline\sigma^l_\gamma \frac{\partial \overline\sigma_\gamma}{\partial x^l} - \sigma^h_\gamma P_k\frac{\partial \sigma^k_\gamma}{\partial x^h} \bigg)}_{\widetilde \mu} \bigg)(Y_t,t) \dif t
	\end{split}
\end{align}
Using \eqref{dP} we can split $\widetilde \mu$ in its orthogonal and tangential components: on $M$ we have
\begin{align}\label{stratdrift}
	\begin{split}
		\widetilde \mu &= \overline \mu + \frac 12 \sum_{\gamma = 1}^n \bigg( \overline \sigma^l_\gamma \bigg( \frac{\partial P_k}{\partial x^l} \sigma^k_\gamma + P_k \frac{\partial \sigma^k_\gamma}{\partial x^l} \bigg) - \sigma^h_\gamma P_k \frac{\partial \sigma^k_\gamma}{\partial x^h} \bigg) \\
		&=  \overline\mu + \frac 12 \sum_{\gamma = 1}^n \bigg( \frac{\partial P_k}{\partial x^l} \overline \sigma^l_\gamma \widecheck \sigma^k_\gamma + \frac{\partial P_k}{\partial x^l} \overline \sigma^l_\gamma \overline \sigma^k_\gamma + \overline \sigma^l_\gamma P_k \frac{\partial \sigma^k_\gamma}{\partial x^l}  - \sigma^h_\gamma P_k \frac{\partial \sigma_\gamma^k}{\partial x^h} \bigg) \\
		&= \underbrace{\overline \mu + \frac 12 \sum_{\gamma = 1}^n \bigg( \frac{\partial^2 \pi}{\partial x^i \partial x^j} \overline \sigma^i_\gamma \widecheck \sigma^j_\gamma - \widecheck \sigma^h_\gamma P_k \frac{\partial \sigma^k_\gamma}{\partial x^h} \bigg)}_{\in TM}  + \underbrace{\frac 12 \sum_{\gamma = 1}^n  \frac{\partial^2 \pi}{\partial x^i \partial x^j}\overline \sigma^i_\gamma \overline \sigma^j_\gamma}_{\in T^\bot M}
	\end{split}
\end{align}
with implied evaluation of all terms at $(y,t)$.

We now move on to the It\^o-jet projection. Let $\bbF \in \Gamma \mathrm{Diff}_\mathrm{Sch}^{\, n} \bbR^d$ as in \eqref{SMembeddedSDE}, so that the Schwartz-Meyer equation it defines coincides with the It\^o equation \eqref{sde}. We can then write \eqref{GeneralSMProj} using matrix notation as
\begin{align}
	\begin{split}
		\begin{bmatrix}
			\dif Y_t \\
			\frac 12 \dif [Y]_t 
		\end{bmatrix} = \begin{bmatrix}
			\frac{\partial \pi}{\partial x} & \frac{\partial^2 \pi}{\partial x^2} \\ 0 &  \frac{\partial \pi}{\partial x} \odot \frac{\partial \pi}{\partial x}
		\end{bmatrix}(Y_t)\begin{bmatrix}
			F & G \\ 0 &  F \odot  F
		\end{bmatrix}(Y_t,t) \begin{bmatrix}
			\dif W_t \\
			\frac 12 \dif [W]_t 
		\end{bmatrix}
	\end{split}	
\end{align}
of which the first line reads
\begin{align}\label{firstline}
	\begin{split}
		\dif Y_t &= \frac{\partial \pi}{\partial x^k}(Y_t) \bigg( F^k_\gamma(Y_t,t)\dif W_t^\gamma + F^k_0(Y_t,t)\dif t + \frac 12 \sum_{\gamma = 1}^n G_{\gamma \gamma}(Y_t,t) \dif t\bigg) \\
		&\mathrel{\hphantom{=}} + \frac 12 \sum_{ \gamma = 1}^n \frac{\partial^2 \pi}{\partial x^i \partial x^j}(Y_t) F^i_\gamma F^j_\gamma(Y_t,t) \dif t \\
		&= \overline \sigma_\gamma(Y_t,t) \dif W^\gamma_t + \bigg( \underbrace{\overline \mu + \frac 12 \sum_{\gamma = 1}^n \frac{\partial^2 \pi}{\partial x^i \partial x^j}\sigma^i_\gamma \sigma^j_\gamma}_{\widehat \mu} \bigg)(Y_t,t) \dif t
	\end{split}
\end{align}
\begin{rem}\label{generator}
	We can write the It\^o-jet-projected drift $\widehat \mu$ as the generator of the SDE, applied to the tubular neighbourhood projection $\pi$: 
	\begin{equation}\label{itojetdrift}
		\widehat \mu(y,t) = \frac{\partial \pi}{\partial x^k} \mu^k(t,y) + \frac 12 \sum_{\gamma = 1}^n \frac{\partial^2 \pi}{\partial x^i \partial x^j}\sigma^i_\gamma \sigma^j_\gamma (t,y) = (\mathscr L_t \pi)(y)
	\end{equation}
\end{rem}
In \cite{AB18} the field of Schwartz morphisms $\bbF$ is taken to be induced by a (time-homogeneous) field of maps $f$ as in \eqref{jets}. In this approach we can use functoriality of $\bbT$ to write
\begin{equation}
	\widehat \bbF(y) = \bbT_y \pi \circ \bbF(y) = \bbT_y \pi \circ \bbT_0 f_{y} =  \bbT_0 (\pi \circ f_y)
\end{equation}
thus obtaining an SDE defined by the field of (2-jets of) maps given by projecting the original field of maps onto $M$ with the tubular neighbourhood projection $\pi$.

Finally, we consider the It\^o-vector projection of \eqref{sde}. By \eqref{itocoefficientsembedded}, in coordinates this amounts to projecting \eqref{sde} to the It\^o SDE on $M$ with diffusion coefficients given by $\overline \sigma_\gamma$ and drift
\begin{equation}\label{generalitovector}
	\overrightarrow \mu = \underbrace{\overline \mu}_{\in TM} + \underbrace{\frac 12 \sum_{\gamma = 1}^n \frac{\partial^2 \pi}{\partial x^i \partial x^j} \overline \sigma^i_\gamma \overline \sigma^j_\gamma}_{\in T^\bot M}
\end{equation}

To summarise, all three projections of the It\^o equation \eqref{sde} agree on how to map the diffusion coefficients, and the orthogonal components of the drift terms will all be fixed by the constraint \eqref{itotan}, while their tangential projections are given by (respectively Stratonovich, It\^o-jet, It\^o-vector)
\begin{equation}\label{mutable} 
	\begin{array}{|l|l|}
		\hline
		P\widetilde \mu & \displaystyle \vphantom{\Bigg[} \overline \mu + \frac 12 \sum_{\gamma = 1}^n \bigg( \frac{\partial^2 \pi}{\partial x^i \partial x^j} \overline \sigma^i_\gamma \widecheck \sigma^j_\gamma - \widecheck \sigma^h_\gamma P_k \frac{\partial \sigma^k_\gamma}{\partial x^h} \bigg)  \\
		\hline
		P\widehat \mu & \displaystyle \vphantom{\Bigg[} \overline \mu + \sum_{\gamma = 1}^n \frac{\partial^2 \pi}{\partial x^i \partial x^j} \overline \sigma^i_\gamma \widecheck \sigma^j_\gamma  \\
		\hline
		P \overrightarrow \mu & \vphantom{\Bigg[} \overline \mu \\
		\hline
	\end{array}
\end{equation}

By calculations similar to \eqref{stratdrift} we can compute the projections of \eqref{StratSDE} in Stratonovich form: again, all three projections will orthogonally project the diffusion coefficients, and behave as follows on the Stratonovich drifts.
\begin{equation}\label{btable} 
	\begin{array}{|l|l|}
		\hline
		\widetilde b & \displaystyle \vphantom{\Bigg[} \overline b  \\
		\hline
		\widehat b & \displaystyle \vphantom{\Bigg[} \overline b + \frac 12 \sum_{\gamma = 1}^n \bigg( \widecheck \sigma^h_\gamma P_k \frac{\partial \sigma^k_\gamma}{\partial x^h} + \frac{\partial^2 \pi}{\partial x^i \partial x^j} \overline \sigma^i_\gamma \widecheck \sigma^j_\gamma \bigg)   \\
		\hline
		\overrightarrow b & \vphantom{\Bigg[} \overline b + \displaystyle\frac 12 \sum_{\gamma = 1}^n \bigg( \widecheck \sigma^h_\gamma P_k \frac{\partial \sigma^k_\gamma}{\partial x^h} - \frac{\partial^2 \pi}{\partial x^i \partial x^j} \overline \sigma^i_\gamma \widecheck \sigma^j_\gamma \bigg) \\
		\hline
	\end{array}
\end{equation}
From now on we will consider \eqref{sde} as being our starting point, unless otherwise mentioned, and thus refer to \eqref{mutable} when comparing the three projections.

We end this section with a brief comparison of the three projections, leaving a detailed analysis of their differences to \autoref{subsec:comparison}. The three projections coincide if $\sigma_\gamma \in TM$ for $\gamma = 1,\ldots, n$ (which includes the ODE case $\sigma_\gamma = 0$), in which case the diffusion coefficients remain unaffected, and the tangent component of the projected drift is simply given by $\overline \mu$. If $\sigma_\gamma \in T^\bot M$ for $\gamma = 1,\ldots, n$ all three projections result in an ODE on $M$, and the It\^o-jet and It\^o-vector projections coincide. Another case in which the It\^o-jet and It\^o-vector projections coincide is when the second derivatives of $\pi$ vanish: this occurs in particular if $M$ is embedded affinely, i.e.\ it coincides with some open set of an affine space of $\bbR^d$. All three projections forget the orthogonal part of the (It\^o or Stratonovich) drift. We observe from \eqref{mutable} that the It\^o-jet and It\^o-vector projections of \eqref{sde} only depend on the values of the It\^o-coefficients on $M$. The Stratonovich projection, instead, could additionally depend on the tangential components of the derivatives of the diffusion coefficients in the direction of their normal components. Naturally, the situation is reversed when projecting \eqref{StratSDE}: here it is the Stratonovich projection that only depends on the values of the coefficients on $M$, while the It\^o-jet and -vector projections might depend on the mentioned derivative term. 
\begin{expl}[The projections in the case $M$ time-dependent]\label{tproj}
	Recalling \autoref{timedepM} (and the map $\widetilde \pi$ defined therein) we may ask whether there is a way to consider the three SDE projections in the case of $M$ time-dependent. The most natural way to define this is to consider, as done in \eqref{tX}, the joint equation satisfied by $(t,X_t)$, project its coefficients in the three ways onto $\widetilde M$, thus obtaining a solution of the form $(t,Y_t)$: this uses that $\widetilde \pi^0(t,y) = t$ (with time the $0^\text{th}$ coordinate), which is instead not necessarily satisfied by the Riemannian tubular neighbourhood projection onto $\widetilde M$. It is easily checked that the formulae \eqref{mutable} for the tangential component of the drift of $Y_t$ continue to hold with the substitution of $\pi_t$ for $\pi$ (so that also the projection onto the tangent space $P$ is now time-dependent), whereas in all three cases the orthogonal component of the drift picks up the term $\dot \pi_t$, needed to keep the process on the evolving manifold $M_t$. In particular, in the It\^o-jet case we have
	\begin{equation}\label{tgenerator}
		\widehat \mu(y,t) = (\mathscr L_t \pi_t)(y) + \dot \pi_t(y) = \widetilde{\mathscr L} \widetilde\pi (t,y)
	\end{equation}
	where $\mathscr L_t$ is the generator of $X$ and $\widetilde{\mathscr L}$ is that of $(t,X_t)$ (which can be considered as being a time-homogeneous Markov process). This identity extends the observation made in \autoref{generator}. The same term $\dot\pi_t$ should be added to the Stratonovich drifts \eqref{btable} for the extension to the case of $M$ time-dependent.
\end{expl}

\section{The optimal projection}\label{sec:optimal}
In the previous section we showed how to abstractly project manifold-valued SDEs onto submanifolds in three (possibly) different ways, and specialised these constructions to the case of $M \hookrightarrow \bbR^d$-valued diffusions. In this section we will seek the \emph{optimal} projection of an SDE for $X_t$, which we write in It\^o form as \eqref{sde}. Let
\begin{equation} \label{Ysde}
	\dif Y_t^k = \accentset{\circ}\sigma^k_\gamma(Y_t,t) \dif W^\gamma_t + \accentset{\circ}\mu^k(Y_t,t) \dif t, \quad X_0 = y_0 \in M
\end{equation}
be the $M$-valued SDE to be defined, which we write in $\bbR^d$-coordinates. Its coefficients $\accentset{\circ}\sigma_\gamma$ and $\accentset{\circ}\mu$ are to be treated as unknowns, to be determined by the optimisation criteria that involve the minimisation of the quantities
\begin{equation}\label{eq:opCriteria}
E[|Y_t - X_t|^2], \qquad E[|Y_t - \pi(X_t)|^2], \qquad |E[Y_t - X_t]|^2
\end{equation}
asymptotically for small $t$. Before we define the optimality criteria precisely, it is important to note that such expectations are undefined if the solution to either SDE is explosive, or, in the second case, even if it exits the tubular neighbourhood of $M$ on which $\pi$ is defined. The problem must be slightly changed so as to ensure that we are minimising a well-defined quantity. One option is to take the above expectations on the event $\{t \leq \tau_r\}$, where 
\begin{equation}\label{tau}
	\tau_r \coloneqq \min \{t \geq 0 : |(X_t,Y_t) - (y_0,y_0)|^2 \geq r^2 \}
\end{equation}
for some suitable $r>0$. However, since for such optimality criteria the values of the vector fields of both SDEs outside the ball $B_{(y_0,y_0)}(r) \subseteq \bbR^{2d}$ are irrelevant, it is simpler to just assume that they vanish outside $B_{(y_0,y_0)}(2r)$. Since the optimisation criteria will only determine the value of $\accentset{\circ}\sigma$, $\accentset{\circ}\mu$ at the initial condition, this is really only an assumption on $\sigma$ and $\mu$. The following proposition reassures us that, at least in well-behaved cases, this does not alter the problem in a way that interferes with the optimisation (which, as will be seen shortly, only involves the Taylor expansions of order 2 of \eqref{eq:opCriteria} in $t = 0$).
\begin{lem}
	Let $X,Y,y_0,\tau_r$ be as above, $U$ a neighbourhood of $(y_0,y_0)$ in $\bbR^{2d}$ and assume that there exists deterministic $\varepsilon>0$ s.t.\ $X_t,Y_t \in U$ for $t \in [0,\varepsilon]$. Let $f \colon U \times [0,\varepsilon] \to \bbR$ be continuous s.t.\ $f(y_0,y_0,0) = 0$, and assume moreover that $E[\max_{0 \leq t \leq \varepsilon}|f(X_t,Y_t,t)|] < \infty$ (this holds, in particular, under the global Lipschitz assumptions that guarantee SDE exactness \cite[Theorem 11.2]{RW00}). Then for any $r > 0$ with $\overline B_r(y_0,y_0) \subseteq U$
	\begin{equation}
		E[f(X_t,Y_t,t)] - E[f(X_t,Y_t,t);t<\tau_r]
	\end{equation}
belongs to $O(t^n)$ for all $n \in \mathbb N$ as $t \to 0$.
\begin{proof}
	Fix $r$, and let $\tau \coloneqq \tau_r$. The It\^o formula yields the a decomposition $|(X_t,Y_t) - (y_0,y_0)|^2 = L_t + A_t$ with $L_t$ sum of Brownian integrals and $A_t$ time integral, all of which for $t \leq \tau \wedge \varepsilon$ have bounded integrand (by continuity of the SDE coefficients and compactness of $\overline B_r(y_0,y_0) \times [0,\varepsilon]$). $[L]_t$ can be expressed as a time integral with bounded integrand: let $R>0$ bound the sum of the absolute values of all integrands mentioned for $t \in [0,\tau \wedge \varepsilon]$. Then, still for $t \leq \tau \wedge \varepsilon$ we have $|A_t|, [L]_t \leq Rt$, and for any $\xi > 0$ it holds that $|(X_t,Y_t) - (y_0,y_0)|^2 \leq L_t + R\xi$ for $0 \leq t \leq \varepsilon \wedge \xi$. Letting $\xi \coloneqq r^2/(3R)$, on $[0,\varepsilon \wedge \xi]$ we have
	\begin{align}\label{taubound}
		\begin{split}
			P[t \geq \tau] &= P\big[\max_{0 \leq s \leq t}|(X_s,Y_s) - (y_0,y_0)|^2 \geq r^2 \big] \\
			&= P\big[\max_{0 \leq s \leq \tau \wedge t}|(X_s,Y_s) - (y_0,y_0)|^2 \geq r^2 \big] \\
			&\leq P\big[\max_{0 \leq s \leq \tau \wedge t} L_s > r^2/2 \big] \\
			&= P\big[\max_{0 \leq s \leq t} L_{\tau \wedge s} > r^2/2 \big] \\
			&\leq \exp\bigg(-\frac{r^4}{4Rt} \bigg)
		\end{split}
	\end{align}
by the tail estimate \cite[Theorem 37.8 p.77]{RW00}. Now, for $t \in [0,\varepsilon \wedge \xi]$ by Cauchy-Schwarz
\begin{equation}
	\begin{split}
		&\big| E[f(X_t,Y_t,t)] - E[f(X_t,Y_t,t);t<\tau] \big| \\
		={} &\big|E[f(X_t,Y_t,t);t\geq\tau]\big| \\
		\leq{}  &E[f(X_t,Y_t,t)^2]^{1/2} P[t \geq \tau]^{1/2} \\
		\leq{}  &E\Big[\max_{[0,t]}f(X_s,Y_s,s)^2\Big]^{1/2} P[t \geq \tau]^{1/2} \\
		\lesssim{} &\exp\bigg(-\frac{r^4}{4Rt} \bigg)
	\end{split}
\end{equation}
since the first factor also vanishes as $t \to 0$, by the hypotheses on $f,X,Y$ and dominated convergence.
\end{proof}
\end{lem}

We proceed with the constrained optimisation problem, assuming all SDE coefficients to be compactly supported; this means all local martingales involved will be martingales, and that we may use Fubini to pass to the expectation inside integrals in $\dif t$. If we can write the Taylor expansion of the strong error
\begin{equation}\label{obvmin}
E \big[ |Y_t - X_t|^2\big] = a_1 t + a_2 t^2 + o(t^2)
\end{equation}
a first goal could be to minimise the leading coefficient $a_1$ (of course there is no constant term because $Y_0 = y_0 = X_0$). Using It\^o's formula, and intending with $\simeq$ equality of differentials up to differentials of martingales, we have
\begin{align*}
		&\hphantom{{}={}}\dif |Y_t - X_t|^2 \\
		&= \dif \sum_{k = 1}^d (Y^k_t - X^k_t)^2 \\
		&= 2 \sum_{k = 1}^d \big( (Y^k_t - X^k_t) \dif Y_t^k - (Y^k_t - X^k_t) \dif X^k_t \big) + \sum_{k = 1}^d \big( \dif Y^k_t \dif Y^k_t + \dif X^k_t \dif X^k_t - 2 \dif X^k_t \dif Y^k_t \big) \\
		&\simeq \sum_{k = 1}^d \bigg[ 2\bigg( \sum_{\gamma = 1}^n \int_0^t \big(\accentset{\circ} \sigma^k_\gamma(Y_s,s) - \sigma^k_\gamma(X_s,s)\big) \dif W^\gamma_s \\
		&\mathrel{\hphantom{\simeq \sum_{\gamma,k} \bigg[ 2\bigg(}}+ \int_0^t \big(\accentset{\circ} \mu^k(Y_s,s) - \mu^k(X_s,s)\big) \dif s \bigg) \big(\accentset{\circ} \mu^k(Y_t,t) - \mu^k(X_t,t)\big)\\
		&\mathrel{\hphantom{\simeq \sum_{\gamma,k} \bigg[ }} + \sum_{\gamma = 1}^n\big( \accentset{\circ} \sigma_\gamma^k (Y_t,t)^2 +  \sigma_\gamma^k(X_t,t)^2 -2 \sigma_\gamma^k(X_t,t) \accentset{\circ} \sigma_\gamma^k(Y_t,t) \big) \bigg] \dif t
\end{align*}
We now compute the expectation:
\begin{align*}
		&\mathrel{\hphantom{=}}E \big[ |Y_t - X_t|^2 \big] \\
		&= \sum_{k = 1}^d 2E \bigg[ \int_0^t \bigg( \sum_{\gamma = 1}^n\int_0^s \big(\accentset{\circ} \sigma^k_\gamma(Y_u,u) - \sigma^k_\gamma(X_u,u)\big) \dif W^\gamma_u \bigg) \big(\accentset{\circ} \mu^k(Y_s,s) - \mu^k(X_s,s)\big) \dif s  \bigg] \\
		&\mathrel{\hphantom{= \sum_{k = 1}^m}} + 2E \bigg[ \int_0^t \bigg( \int_0^s \big( \accentset{\circ} \mu^k(Y_u,u) - \mu^k(X_u,u) \big) \dif u \bigg)  \big( \accentset{\circ} \mu^k(Y_s,s) - \mu^k(X_s,s) \big) \dif s\bigg] \\
		&\mathrel{\hphantom{= \sum_{k = 1}^m}} + E \bigg[ \sum_{\gamma = 1}^n\int_0^t \big( \accentset{\circ} \sigma_\gamma^k \accentset{\circ} \sigma_\gamma^k(Y_s,s) +  \sigma_\gamma^k  \sigma_\gamma^k(X_s,s) -2 \sigma_\gamma^k(X_s,s) \accentset{\circ} \sigma_\gamma^k(Y_s,s) \big) \dif s \bigg] \\
		&=  \int_0^t E\bigg[ \sum_{k = 1}^d 2\bigg( \sum_{\gamma = 1}^n\int_0^s \big(\accentset{\circ} \sigma^k_\gamma(Y_u,u) - \sigma^k_\gamma(X_u,u)\big) \dif W^\gamma_u \bigg) \big(\accentset{\circ} \mu^k(Y_s,s) - \mu^k(X_s,s)\big) \\
		&\mathrel{\hphantom{=  \int_0^t E\bigg[ \sum_{\gamma,k}}} + 2\bigg( \int_0^s \big( \accentset{\circ} \mu^k(Y_u,u) - \mu^k(X_u,u) \big) \dif u \bigg)  \big( \accentset{\circ} \mu^k(Y_s,s) - \mu^k(X_s,s) \big) \\
		&\mathrel{\hphantom{=  \int_0^t E\bigg[ \sum_{\gamma,k}}} +\sum_{\gamma = 1}^n  \big( \accentset{\circ} \sigma^k_\gamma(Y_s,s) - \sigma^k_\gamma(X_s,s) \big)^2
		\bigg] \dif s \\
		&\eqqcolon \int_0^t E[Z_s] \dif s
\end{align*}
and differentiating, with reference to \eqref{obvmin} we have 
\begin{equation}\label{a1}
	a_1 = \frac{\dif}{\dif t} \bigg|_0^+ \int_0^t E[Z_s] \dif s = \sum_{\gamma = 1}^n \abs{ \accentset{\circ} \sigma_\gamma(y_0,0) - \sigma_\gamma(y_0,0) }^2
\end{equation}
Since $a_1$ only depends on the diffusion coefficients, its minimisation is expressed by the constrained optimisation problem whose solution is simply given by projecting the $\sigma_\gamma$'s onto $TM$:
\begin{align}\label{constraeasy}
	\begin{dcases}
		\text{minimise } \sum_{\gamma = 1}^n \abs{ \accentset{\circ} \sigma_\gamma - \sigma_\gamma }^2 \\
		\text{subject to }
		Q^k_h\accentset{\circ} \sigma^h_\gamma = 0 
	\end{dcases} \Longleftrightarrow \accentset{\circ} \sigma = \overline \sigma = P\sigma
\end{align}
Here we have omitted evaluation at the initial condition $(0,y_0)$. Since we have not obtained a condition on $\accentset{\circ} \mu$ our SDE \eqref{Ysde} is still underdetermined, and the condition would be satisfied by the Stratonovich projection of \eqref{sde}.

One idea to obtain a condition on $\accentset{\circ} \mu$ would be to minimise $a_2$ in \eqref{obvmin}. This attempt, however, has the drawback that we are minimising the second Taylor coefficient of a function without its first vanishing (unless the $\sigma_\gamma$'s are already tangent to start with: in this case the minimisation of $a_2$ can be seen to result in the three projections, which all coincide). Although this approach is discussed in \cite{ABR19}, we will not do so here, as there are more sound optimisation criteria. Indeed, we can look at the Taylor expansion of the weak error
\begin{equation}\label{weakerror}
	\abs{E[Y_t - X_t ]}^2 = b_2 t^2 + o(t^2) \quad \text{as } t \to 0^+
\end{equation}
We compute the term on the left as
\begin{align}
	\begin{split}
		\abs{E[Y_t - X_t]}^2 &= \bigg\lvert \int_0^t E[ \accentset{\circ} \mu(Y_s,s) - \mu(X_s,s) ] \dif s \bigg\rvert^2
	\end{split}
\end{align}
and
\begin{align}
	\begin{split}
		&\mathrel{\hphantom{=}}\frac{\dif}{\dif t} \bigg\lvert \int_0^t E[ \accentset{\circ} \mu(Y_s,s) - \mu(X_s,s) ] \dif s \bigg\rvert^2 \\
		&= 2 E[ \accentset{\circ} \mu(Y_t,t) - \mu(X_t,t) ]\int_0^t E[ \accentset{\circ} \mu(Y_s,s) - \mu(X_s,s) ] \dif s \\
		&\mathrel{\hphantom{=}}\frac{\dif^2}{\dif t^2}\bigg|_0 \bigg\lvert \int_0^t E[ \accentset{\circ} \mu(Y_s,s) - \mu(X_s,s) ] \dif s \bigg\rvert^2 = 2\abs{\accentset{\circ} \mu(y_0,0) - \mu(y_0,0)}^2
	\end{split}
\end{align}
which confirms that \eqref{weakerror} lacks a linear term, and we have
\begin{equation}
	b_2 = \abs{\accentset{\circ} \mu - \mu}^2
\end{equation}
Requiring the minimisation of $b_2$ is thus independent of the minimisation of $a_1$ above, and results in the constrained optimisation problem
\begin{align}\label{constravector}
	\begin{dcases}
		\text{minimise } \abs{\accentset{\circ} \mu-\mu}^2 \\
		\text{subject to }
		Q^k_h\accentset{\circ} \mu^h_\gamma = \frac 12 \sum_{\gamma = 1}^n \frac{\partial^2 \pi^k}{\partial x^i \partial x^j} \overline\sigma^i_\gamma \overline\sigma^j_\gamma 
	\end{dcases} \Longleftrightarrow \accentset{\circ} \mu = \overline \mu + \frac 12 \sum_{\gamma = 1}^n \frac{\partial^2 \pi}{\partial x^i \partial x^j} \overline\sigma^i_\gamma \overline\sigma^j_\gamma
\end{align}
A quick glance at \eqref{mutable} shows we have proven the following
\begin{thm}[Optimality of the It\^o-vector projection]\label{maintheoremV}
	The coefficients $\accentset{\circ}\sigma_\gamma, \accentset{\circ}\mu$ of the $M$-valued SDE \eqref{Ysde} that solve the constrained optimisation problem
	\begin{equation}\label{constrthmV}
		\begin{dcases}
			\emph{minimise $a_1$ in  \eqref{obvmin} and $b_2$ in \eqref{weakerror}} \\
			\emph{subject to \eqref{itotan}} 
		\end{dcases}
	\end{equation}
	for all initial conditions $X_0 = Y_0 = y_0 \in M$ are given (uniquely for $t=0$) by the It\^o-vector projection of the original SDE \eqref{sde}.
\end{thm}

\begin{rem}\label{remtdep}
	In defining the three projections in \autoref{sec:projecting} we intended for the projected coefficients to still be time-dependent if the original ones were. The optimality requirement only fixes the coefficients at the initial condition, at time $0$, i.e.\ $\accentset{\circ} \sigma_\gamma(y_0,0), \accentset{\circ} \mu(y_0,0)$. To retain the time-dependence we may consider the optimisation involving all time-translated initial conditions $Y_{t_0} = y_0$.
\end{rem}

\begin{rem}\label{startfromstrat}
	Note that the form (It\^o or Stratonovich) the initial SDE is provided in is irrelevant: if we had begun with \eqref{StratSDE} instead of \eqref{sde} the optimality criterion would still have led us to the It\^o-vector projection, which for the Stratonovich drift would have taken the form $\overrightarrow b$ in \eqref{btable}. The only reason to start with an It\^o SDE is that the calculations are simpler, and it is possible to express the optimal coefficients as functions of the values of the coefficients of the original SDE, without reference to their derivatives.
\end{rem}

$ $\par
The optimisation of \autoref{maintheoremV} has the disadvantage of coming from the two separate minimisations of $a_1$ and $b_2$, which are Taylor coefficients of different quantities. There is a different way of arriving at coefficients by successively minimising the Taylor coefficients of the same quantity, with the first minimisation resulting in a null term. The idea is to consider
\begin{equation} \label{projmin}
	E\big[|Y_t - \pi(X_t)|^2 \big] = c_1 t + c_2 t^2 + o(t^2)
\end{equation}
where $X,Y,\tau$ are respectively as in \eqref{sde}, \eqref{Ysde},\eqref{tau}, with the requirement on $r$ that $B_r(y_0)$ be contained in the domain of $\pi$. The map $\pi$ is the one defined in \eqref{defpi}, although it can more generally satisfy \eqref{pigeneral}. Letting $\accentset{\circ} \sigma_\gamma, \accentset{\circ}\mu$ resume their status as unknowns, we proceed with the calculations.
\begin{align*}
		&\hphantom{{}={}}\dif |Y_t - \pi(X_t)|^2 \\
		&= \dif \sum_{k = 1}^d (Y^k_t - \pi^k(X_t))^2 \\
		&=  \sum_{k = 1}^d \bigg[2(Y^k_t - \pi(X^k_t)) \dif Y_t^k - 2(Y^k_t - \pi(X^k_t)) \frac{\partial \pi^k}{\partial x^h}(X_t) \dif X^h_t + \dif Y^k_t \dif Y^k_t \\
		&\mathrel{\hphantom{=\sum_{k = 1}^d}} + \bigg( \frac{\partial \pi^k}{\partial x^i}\frac{\partial \pi^k}{\partial x^j}(X_t) - (Y^k_t - \pi^k(X_t))\frac{\partial^2\pi^k}{\partial x^i\partial x^j}(X_t) \bigg) \dif X^i_t \dif X^j_t - 2 \frac{\partial \pi^k}{\partial x^h}(X_t) \dif X^h_t \dif Y^k_t \bigg] \\
		&\simeq \sum_{k = 1}^d \bigg[ 2\big(Y^k_t - \pi^k(X_t)\big) \bigg(\accentset{\circ} \mu^k(t,Y_t) - \frac{\partial \pi^k}{\partial x^h}\mu^h(t,X_t) - \frac 12 \sum_{\gamma = 1}^n \frac{\partial^2 \pi^k}{\partial x^i \partial x^j}\sigma_\gamma^i \sigma_\gamma^j(X_t,t) \bigg)\\
		&\mathrel{\hphantom{\simeq \sum_{k = 1}^d \bigg[ }} + \sum_{\gamma = 1}^n \bigg( \accentset{\circ} \sigma_\gamma^k \accentset{\circ} \sigma_\gamma^k(t,Y_t) + \frac{\partial \pi^k}{\partial x^i} \frac{\partial \pi^k}{\partial x^j} \sigma_\gamma^i \sigma_\gamma^j(t,X_t) -2 \frac{\partial \pi^k}{\partial x^h}\sigma_\gamma^h(t,X_t) \accentset{\circ} \sigma_\gamma^k(t,Y_t) \bigg) \bigg] \dif t \\
		& \eqqcolon Z_t \dif t
\end{align*}
and
\begin{align}\label{somelabel}
	\begin{split}
		&\mathrel{\hphantom{\approx}}\frac{\dif}{\dif t} E\bigg[ \int_0^t Z_s \dif s \bigg]	\\
		&= E \bigg[ \sum_{k = 1}^d 2\big(Y^k_t - \pi^k(X_t)\big) \bigg(\accentset{\circ} \mu^k(t,Y_t) - \frac{\partial \pi^k}{\partial x^h}\mu^h(t,X_t) - \frac 12 \sum_{\gamma = 1}^n  \frac{\partial^2 \pi}{\partial x^i \partial x^j}\sigma_\gamma^i \sigma_\gamma^j(X_t,t) \bigg) \\
		&\mathrel{\hphantom{\simeq E \bigg[ \sum_{k = 1}^d }} +\sum_{\gamma = 1}^n \bigg( \accentset{\circ} \sigma_\gamma^k(Y_t,t) - \frac{\partial \pi^k}{\partial x^h} \sigma_\gamma^h(X_t,t) \bigg)^2 \bigg]
	\end{split}
\end{align}
and therefore
\begin{align}
	\begin{split}
		c_1 = \frac{\dif}{\dif t}\bigg|_0^+ E\bigg[ \int_0^t Z_s \dif s \bigg] = \sum_{\gamma = 1}^n |\accentset{\circ} \sigma_\gamma - P\sigma_\gamma |^2
	\end{split}
\end{align}
(evaluation at $(y_0,0)$ is implied). Thus $c_1$ vanishes if and only if $\accentset{\circ} \sigma \coloneqq P\sigma$. Continuing as before and we have
\begin{align}
	\begin{split}
		\dif Z_t &\simeq  \sum_{k = 1}^d 2 \big(Y^k_t - \pi^k(X_t)\big) \dif (\ldots) \\
		&\mathrel{\hphantom{\simeq  \sum_{\gamma, k} }}+ 2 \bigg(\accentset{\circ} \mu^k(t,Y_t) - \frac{\partial \pi^k}{\partial x^h}\mu^h(t,X_t) - \frac 12 \frac{\partial^2 \pi^k}{\partial x^i \partial x^j}\sigma_\gamma^i \sigma_\gamma^j(X_t,t) \bigg)^2 \dif t \\
		&\mathrel{\hphantom{\simeq  \sum_{\gamma, k} }} + 2\sum_{\gamma = 1}^n \bigg( \accentset{\circ} \sigma^k_\gamma(Y_t,t) - \frac{\partial \pi^k}{\partial x^h} \sigma_\gamma^h(X_t,t) \bigg) \dif (\ldots) \\
		&\mathrel{\hphantom{\simeq}} +  2 f( \sigma, J \sigma, H \sigma;\accentset{\circ} \sigma, J \accentset{\circ} \sigma, H \accentset{\circ} \sigma; \mu, J \mu)|_{X_t,Y_t,t} \dif t 
	\end{split}
\end{align}
for some smooth function $f$ ($J$ denotes Jacobian and $H$ Hessian), which we denote $f_t$ for short; the differentials $\dif (\ldots)$ can be ignored, since their factors will vanish when evaluated below.
\begin{align}
	\begin{split}
		c_2 &= \frac 12 \frac{\dif^2}{\dif t^2}\bigg|_0^+ E\big[Z_t \big] \\
		&= \sum_{\gamma, k} \bigg(\accentset{\circ} \mu^k - \frac{\partial \pi}{\partial x^h}\mu^h - \frac 12 \frac{\partial^2 \pi}{\partial x^i \partial x^j} \sigma^i_\gamma \sigma^j_\gamma \bigg)^2 + f_t
	\end{split}
\end{align}
The constrained optimisation problem for the minimisation of $c_2$ conditional on the previous minimisation of $c_1$ is thus given by
\begin{align}\label{itojetproblem}
	\begin{split}
		&\begin{dcases}
			\text{minimise }\sum_{k = 1}^d \bigg(\accentset{\circ} \mu^k - \frac{\partial \pi^k}{\partial x^h}\mu^h - \frac 12 \sum_{\gamma = 1}^n \frac{\partial^2 \pi^k}{\partial x^i \partial x^j} \sigma^i_\gamma \sigma^j_\gamma \bigg)^2 \\
			\text{subject to } Q^k_h \accentset{\circ} \mu^h - \frac 12 \sum_{\gamma = 1}^n \frac{\partial^2 \pi^k}{\partial x^i \partial x^j} \accentset{\circ} \sigma^i_\gamma \accentset{\circ} \sigma^j_\gamma = 0
		\end{dcases}	\\
		&\begin{dcases}
			2\bigg(\accentset{\circ} \mu^h - \frac{\partial \pi^h}{\partial x^l}\mu^l - \frac 12 \sum_{\gamma = 1}^n \frac{\partial^2 \pi^h}{\partial x^i \partial x^j} \sigma^i_\gamma \sigma^j_\gamma \bigg) - \sum_{k = 1}^d Q^k_h \lambda^k = 0 \\Q^k_h \accentset{\circ} \mu^h - \frac 12 \sum_{\gamma = 1}^n \frac{\partial^2 \pi^k}{\partial x^i \partial x^j} \accentset{\circ} \sigma^i_\gamma \accentset{\circ} \sigma^j_\gamma = 0
		\end{dcases}\\
		\lambda &\in T_yM, \quad  \mu = P \accentset{\circ} \mu + \frac 12 \sum_{\gamma = 1}^n  \frac{\partial^2 \pi}{\partial x^i \partial x^j} \sigma_\gamma^i \sigma_\gamma^j
	\end{split}
\end{align}
Comparing with \eqref{itojetdrift} we see that we have proven the following
\begin{thm}[Optimality of the It\^o-jet projection]\label{maintheorem}
	The coefficients $\accentset{\circ}\sigma_\gamma, \accentset{\circ}\mu$ of the $M$-valued SDE \eqref{Ysde} that solve the constrained optimisation problem
	\begin{equation}\label{constrthmV}
		\begin{dcases}
			\emph{minimise $c_1$ and $c_2$, conditionally on the minimisation of $c_1$, in \eqref{projmin}} \\
			\emph{subject to \eqref{itotan}} 
		\end{dcases}
	\end{equation}
	for all initial conditions $X_0 = Y_0 = y_0 \in M$ are given (uniquely for $t=0$) by the It\^o-jet projection of the original SDE \eqref{sde}.
\end{thm}
Remarks analogous to \autoref{remtdep} and \autoref{startfromstrat} hold for \autoref{maintheorem}. The It\^o-vector and It\^o-jet projection therefore satisfy different optimality properties, while the Stratonovich projection is suboptimal in both senses. We end the section with the extension of the optimisations to the case of $M$ time-dependent.
\begin{expl}[Optimality for $M$ time-dependent]\label{optMt} Recall the case in which the submanifold $M$ depends smoothly on time, for which we can define similar versions of all three projections \autoref{tproj}. For \autoref{maintheoremV} the optimisation criterion does not require reformulation, while the constraint is modified as described in \autoref{timedepM}: therefore the It\^o-vector projection remains optimal in the case of $M$ time-dependent. For \autoref{maintheoremV} the natural generalisation is given by substituting $\pi_t$ for $\pi$ in \eqref{projmin}. Since $|y - \pi_t(x)| = |(t,y) - \widetilde \pi(t,x)|$, by the definition of the It\^o-jet projection in the case of $M$ time-dependent (and since the calculations in this section never relied on $\pi$ being the Riemannian tubular neighbourhood projection), we have that the time-dependent It\^o-jet projection \eqref{tgenerator} is optimal in this case too.
\end{expl}

\section{Further considerations}\label{sec:further}
In this final section we dig deeper into the details surrounding the It\^o and Stratonovich projections of SDEs, and answer a few lingering questions.

\subsection{Differences between the projections}\label{subsec:comparison}
In this subsection we will provide examples to justify our claim that the vertical rectangles of \eqref{NonCommutativePrism} do not commute.

We begin with an example in which the It\^o-jet and -vector projections coincide, but are different from the Stratonovich projection. This example also shows how the dependence of the Stratonovich projection of \eqref{sde} on the derivatives of the diffusion coefficients can be non-trivial.
\begin{expl} \label{RR2}
	Take $M = \{(x,0) : x \in \bbR\} \hookrightarrow \bbR^2$, $n = 1$ and the It\^o SDEs
	\begin{align}
		\dif \begin{bmatrix}
			X_t \\ Y_t
		\end{bmatrix} = \begin{bmatrix}
			Y_t \\ X_t
		\end{bmatrix} \dif W_t, \quad \dif \begin{bmatrix}
			X_t \\ Y_t
		\end{bmatrix} = \begin{bmatrix}
			0 \\ X_t
		\end{bmatrix} \dif W_t
	\end{align}
	whose diffusion coefficients coincide, and are orthogonal to $M$, on $M$. Their Stratonovich projections onto the affine subspace $M = \bbR$ are respectively given by the ODEs
	\begin{equation}
		\dot{X_t} = - \frac 12 X_t, \quad \dot{X_t} = 0
	\end{equation}
	The It\^o-jet and -vector projections of the two SDEs above coincide (since their coefficients on $M$ coincide) and are trivial. An example where It\^o-jet = It\^o-vector $\neq$ Stratonovich, and where the It\^o projections are non-trivial can be obtained from this by increasing $n$ to $2$ and adding a tangent diffusion coefficient.
\end{expl}

Next, we ask the question of when the Stratonovich and It\^o-jet projections coincide. The following criterion is a rephrasing of \cite[Theorem 5.1]{ABR19}.
\begin{rem}[Fibering property]\label{fibering}
	In general the difference of the Stratonovich- and It\^o-jet-projected drift can be written as
	\begin{equation}\label{ItominusStrat}
		TM \ni \widehat \mu - \widetilde \mu = \frac 12 \sum_{\gamma = 1}^n \bigg(\frac{\partial^2 \pi}{\partial x^i \partial x^j} \overline \sigma^i_\gamma \widecheck \sigma^j_\gamma + \widecheck \sigma^k_\gamma \frac{\partial \pi}{\partial x^h} \frac{\partial \sigma^h_\gamma}{\partial x^k} \bigg)= \frac 12 \sum_{\gamma = 1}^n \frac{\partial}{\partial x^k} \bigg( \frac{\partial \pi}{\partial x^h} \sigma_\gamma^h \bigg)\widecheck \sigma^k_\gamma
	\end{equation}
	Therefore, if we assume that	
	\begin{equation}\label{sigmafibering}
		\frac{\partial \pi}{\partial x^h}(x) \sigma_\gamma^h(x,t) \quad \text{is independent of } x \in \pi^{-1}(\pi(x))
	\end{equation}
	for $x$ in a neighbourhood of $M$ (again, if we are only interested in starting our equation at time zero, the above requirement only needs to be considered for $t = 0$), the derivative of the above quantity along any vector tangent to the fibre of $\pi$ (which at points in $M$ means orthogonal to $M$) vanishes: this means \eqref{ItominusStrat} vanishes and the It\^o jet and Stratonovich projections are equal. Moreover, if, representing the original SDE in Stratonovich form as \eqref{StratSDE}, we additionally have that
	\begin{equation}\label{bfibering}
		\frac{\partial \pi}{\partial x^k}(x) b^k(x,t) \quad \text{is independent of } x \in \pi^{-1}(\pi(x))
	\end{equation}
	then it is immediate to verify that $\pi(X_t)$ is a solution of the Stratonovich projection, and therefore that, letting $Y$ be the solution to the Stratonovich=It\^o-jet projection
	\begin{equation}\label{ypix}
		Y_t = \pi(X_t)
	\end{equation}
	up to the exit time of $X_t$ from the tubular neighbourhood in which $\pi$ is defined. Observe that in the absence of these conditions we cannot expect, in general, to obtain a closed form SDE for $\pi(X_t)$, as the coefficients will depend explicitly on $X_t$. This is even true if \eqref{sigmafibering} holds but \eqref{bfibering} does not, as can be shown simply by considering the ODE case $\sigma_\gamma = 0$.
\end{rem}

\begin{expl}\label{crossdiffusion}
	Let $M = \{(x,y) \in \bbR^2 : x^2 + y^2 = 1\} \hookrightarrow \bbR^2$. $\pi$ is defined in $\bbR^2 \setminus \{0\}$ as $\pi(x,y) = (x^2 + y^2)^{-1/2}(x,y)$. Consider the SDE, dependent on the real parameter $a$
	\begin{align}\label{crossdiff}
		\dif \begin{bmatrix}
			X_t \\ Y_t
		\end{bmatrix} = (X_t^2 + Y_t^2)^a\begin{bmatrix}
			Y_t \\ X_t
		\end{bmatrix} \dif W_t
	\end{align}
	There is a single diffusion coefficient $\sigma$, decomposed as
	\begin{equation}
		\sigma(x,y) = \underbrace{(x^2 + y^2)^{a-1}  (x^2-y^2)\begin{bmatrix}
				-y \\ x
		\end{bmatrix}}_{\overline\sigma(x,y) \in T_{(x,y)}M} + \underbrace{(x^2 + y^2)^{a-1}  2xy \begin{bmatrix}
				x \\ y
		\end{bmatrix}}_{\widecheck \sigma(x,y) \in T^\bot_{(x,y)}M}
	\end{equation}
	Moreover, for $(x,y) \in M$ we have
	\begin{equation}
		J \sigma(x,y) = 	\begin{bmatrix}
			2axy & 2ay^2+1 \\ 2ax^2+1 & 2axy
		\end{bmatrix}
	\end{equation}
	We have
	\begin{align}
		\begin{split}
			J\pi(x,y) &= (x^2+y^2)^{-3/2} \begin{bmatrix} y^2 & -xy \\
				-xy & x^2 \end{bmatrix} \\
			H \pi^1(x,y) &= (x^2+y^2)^{-5/2}\begin{bmatrix}
				-3xy^2 & 2x^2y-y^3 \\ 2x^2y-y^3 & 2xy^2-x^3
			\end{bmatrix} \\
			H \pi^2(x,y) &= (x^2+y^2)^{-5/2}\begin{bmatrix}
				2x^2y-y^3 & 2xy^2-x^3 \\ 2xy^2 - x^3 & -3x^2y
			\end{bmatrix}
		\end{split}
	\end{align}
	We compute, for $(x,y) \in M$
	\begin{align}\label{wholeAfternoon}
		\begin{split}
			\frac{\partial^2 \pi}{\partial x^i \partial x^j} \overline \sigma^i \widecheck \sigma^j  &= -2xy(x^2-y^2)\begin{bmatrix}
				-y \\ x
			\end{bmatrix} \\
			\widecheck \sigma^h P_k \frac{\partial \sigma^k}{\partial x^h} &= 2xy(2ax^4-2ay^4+x^2-y^2)\begin{bmatrix}
				-y \\ x
			\end{bmatrix}
		\end{split}
	\end{align}
	We examine more closely the cases $a = 0$, $a = -1$ and $a = 1$. In the first case (already examined in \cite[\S 5]{ABR19}), the two terms of \eqref{wholeAfternoon} sum to zero, so that \eqref{ItominusStrat} vanishes and the Stratonovich and It\^o-jet projections coincide. Indeed, the fibering property of \eqref{sigmafibering} is verified, as it is easy to see $J\pi(\lambda x, \lambda y) \sigma(\lambda x, \lambda y)$ does not depend on $\lambda >0$. Moreover, since the Stratonovich drift of the equation is given by $-\frac 12 (x,y)$ also \eqref{bfibering} holds and the solution to the Stratonovich=It\^o-jet-projected SDE equals the projection of the solution of the original SDE up to the (a.s.\ infinite) time it hits the origin. However the It\^o-vector projection is distinct, which can be seen by observing that $P\widetilde \mu = P \widehat \mu$ (given by the first term in \eqref{wholeAfternoon}) does not vanish, e.g.\ at the point $(\cos(\pi/6), \sin(\pi/6))$. If $a = -1$ the two terms in \eqref{wholeAfternoon} coincide on $M$ and therefore the Stratonovich projection is identical to the It\^o-vector projection. The It\^o-jet projection, however, is different, again by the nonvanishing of the first term in \eqref{wholeAfternoon} at $(\cos(\pi/6), \sin(\pi/6))$. To generate a case where all three projections are distinct take $a = 1$: all identities can be seen not to hold at the point $(\cos(\pi/6), \sin(\pi/6))$. This case shows that the only projection that preserves the local martingale property is the It\^o-vector.
\end{expl}

\begin{expl}\label{sigmat}
	Consider the case in which $\sigma_\gamma(x,t) = \sigma_\gamma(t)$ do not depend on the state of the solution. In this case, even if \eqref{sde} and \eqref{StratSDE} are equivalent, the projections may still be all different. \eqref{mutable} however shows that 
	\begin{equation}
		\widehat\mu - \overrightarrow \mu = 2(\widetilde \mu - \overrightarrow \mu)
	\end{equation}
	so that if any two projections coincide, they must all.
	An example where all projections are different is given by taking $M$, $d$ as in \autoref{crossdiffusion} and the single, constant diffusion coefficient $\sigma = (1,1)$: all projections differ, for instance at the point $(1,0)$. An example where the projections all coincide is when $n = d$ and $\sigma^k_\gamma = \kron^k_\gamma$:
	\begin{equation}
		\textstyle \sum_\gamma \overline \sigma^i_\gamma \widecheck \sigma^j_\gamma = \sum_\gamma  P^i_\alpha \kron^\alpha_\gamma Q^j_\beta \kron^\beta_\gamma = \sum_\gamma  P^i_\gamma Q^j_\gamma = P^i_\gamma Q^\gamma_j = 0
	\end{equation}
	If the original drift also vanishes, we are in the presence of the trivial SDE for Brownian motion, whose It\^o and Stratonovich projections coincide with the process $\pi(W_t)$ up to the exit time of $W$ from the domain of $\pi$, by the same reasoning of \autoref{fibering}.
\end{expl}

\begin{figure}[ht!]
	\minipage{0.5\textwidth}
	\includegraphics[width=\linewidth]{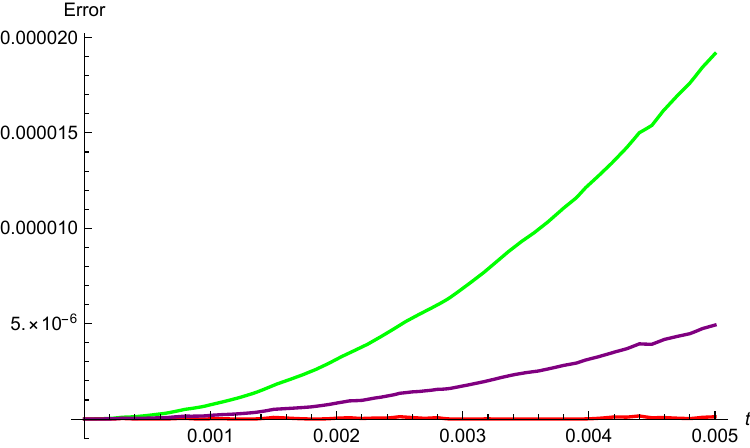}
	\vspace{-20pt}
	\caption*{Weak error \eqref{weakerror}}
	\endminipage\hfill
	\minipage{0.5\textwidth}
	\includegraphics[width=\linewidth]{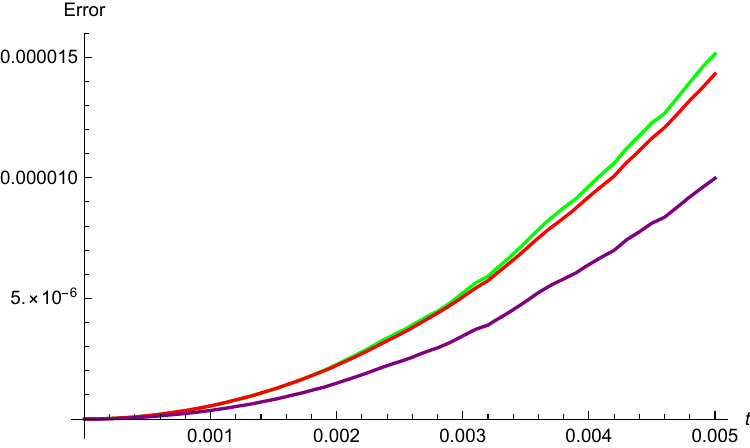}
	\vspace{-22pt}
	\caption*{Mean-squared error \eqref{projmin}}
	\endminipage
	\\
	\minipage{0.85\textwidth}
	\includegraphics[width=\linewidth]{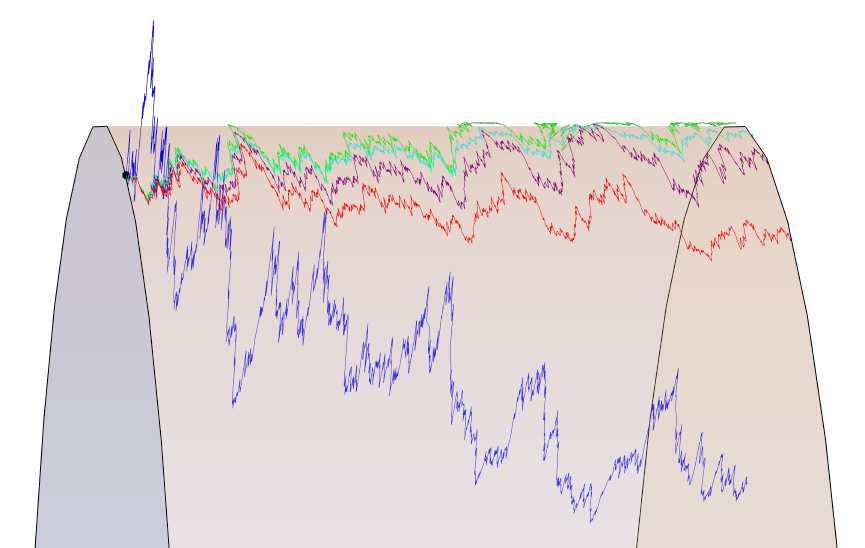}
	\vspace{-20pt}
	\endminipage\hfill
	\minipage{0.15\textwidth}
	\includegraphics[width=\linewidth]{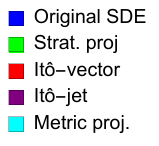}
	\vspace{-20pt}
	\endminipage
	\caption{In these figures we focus on [\autoref{crossdiffusion}, $a = 1$], with initial condition $(\cos(\pi/6),\sin(\pi/6))$, so that all projections are distinct. The two graphs above are respectively plots of the errors $|E[Y_t - X_t]|^2$ and $E[|Y_t - \pi(X_t)|^2]$ for $Y_t$ the solution to the Stratonovich, It\^o-vector and It\^o-jet projections, with the expectation taken over $10^4$ sample paths. We see confirmation of the fact that the It\^o-vector projection performs better in the first error metric, that the It\^o-jet projection does so in the second, and that the Stratonovich projection is markedly suboptimal in both senses (especially in the first, while in the second case it performs very similarly to the It\^o-vector projection). The analogous plot for the error \eqref{obvmin} is not included, as the results for the three projections are visually indistinguishable, in accordance with the fact that all three projections minimise $a_1$ (without it vanishing in this case). The figure below displays one sample path $(t,Y_t)$ where $Y_t$ is each of the following processes: the solution to the original SDE, to the three projected SDEs, and the metric projection $\pi$ applied to the original solution. All sample paths are derived from the same random seed. Since the optimality criteria all involve taking expectation, we do not expect to be able to derive meaningful intuition from a single path, but it is nonetheless informative to have visual confirmation that all projections are distinct, but related.}\label{fig:projComp}
\end{figure}

\pagebreak
In this section we have developed examples that cover all possible situations involving identities, and lack thereof, between the three projections. We summarise them in the table below:
\begin{center} 
	\begin{tabular}[tp]{|l|l|}
		\hline
		$\widetilde \mu = \widehat \mu = \overrightarrow \mu$ &  $\sigma_\gamma \in TM$ and [\autoref{sigmat}, second case] \\
		\hline
		$\widetilde \mu \neq \widehat \mu = \overrightarrow \mu$ &  
		[\autoref{RR2}, first SDE] \\
		\hline
		$\widetilde \mu = \widehat \mu \neq \overrightarrow \mu$ & [\autoref{crossdiffusion}, $a = 0$]  \\
		\hline
		$\widetilde \mu = \overrightarrow \mu \neq \widehat \mu $ & [\autoref{crossdiffusion}, $a = -1$] \\
		\hline
		$\widetilde \mu \neq \widehat \mu \neq \overrightarrow \mu \neq \widetilde \mu$ & [\autoref{crossdiffusion}, $a = 1$] and [\autoref{sigmat}, first case] \\
		\hline
	\end{tabular}
\end{center}

\subsection{Intrinsic optimality of the It\^o projections}\label{subsec:geodesic}

The fact that in \eqref{projmin} we are comparing two points, $Y_t$ and $\pi(X_t)$, which lie in $M$ opens up the possibility of substituting the Euclidean distance with the Riemannian distance of $M$, $\mathscr d_M$, inside the expectation. One can then ask whether this leads to a different optimisation. Let $U$ be a neighbourhood of the initial condition $y_0$ in $\bbR^d$, $V \coloneqq U \cap M$, $\varphi \colon V \to \bbR^m$ a normal chart centred in $y_0$, $\overline \varphi \coloneqq \varphi \circ \pi \colon U \to \bbR^m$. This means that if $G_t$ is a geodesic in $M$ starting at $y_0$, $\varphi (G_t) = vt$ where $\bbR^m \ni v = T_{y_0} \varphi (\dot G_0)$. As a consequence we have that, if $W_{y_0} \in T_{y_0} M$, picking the geodesic $G$ with $G_0 = y_0$, $\dot G_{y_0} = W_{y_0}$, we have that
\begin{equation}\label{d2phi}
	0 = \frac{\dif^2}{\dif t^2}\bigg|_0 \varphi(G_t) = \frac{\dif^2}{\dif t^2}\bigg|_0 \overline\varphi(G_t) = \frac{\partial^2 \overline \varphi}{\partial x^i \partial x^j}(y_0) \dot G_0^i \dot G_0^j + \frac{\partial \overline \varphi}{\partial x^k}(y_0)\ddot G_0^k = \frac{\partial^2 \overline \varphi}{\partial x^i \partial x^j}(y_0) W_{y_0}^i W_{y_0}^j
\end{equation}
since the acceleration of $G$ is orthogonal to $M$. Now, the problem consists of choosing $\accentset{\circ} \sigma_\gamma$ and $\accentset{\circ} \mu$ in such a way that $c'_1$ vanishes and $c'_2$ is minimal in 
\begin{equation} \label{dMmin}
	E\big[{^\varphi \! \!}\mathscr d_M \big(\varphi(Y_t),\varphi(\pi(X_t)) \big)^2  \big] = c'_1 t + c'_2 t^2 + o(t^2)
\end{equation}
where ${^\varphi \! \!}\mathscr d_M (a,b) \coloneqq \mathscr d_M(\varphi^{-1}(a), \varphi^{-1}(b))$. We have expressed $\mathscr d_M$ in normal coordinates in order to be able to use the estimates of \cite[Appendix A]{N14}, which tell us that the derivatives of orders $\leq 3$ of ${^\varphi \! \!}\mathscr d_M$ agree with those of the squared distance function of $\bbR^m$ (in particular those of order 1 and 3 vanish). Since we are only interested in $c'_1$ and $c'_2$, this means we can substitute the LHS of \eqref{dMmin} with
\begin{equation} \label{newdMmin}
	E\big[\mathscr |\varphi(Y_t),\varphi(\pi(X_t))|^2 \big]
\end{equation}
Proceeding as in the computations of \autoref{sec:optimal}, we see that
\begin{equation}
	c_1 = \sum_{\gamma = 1}^n |J\varphi \accentset{\circ} \sigma_\gamma - J\varphi P\sigma_\gamma |^2
\end{equation} 
This quantity is made to vanish exactly as before, namely in the unique case $\accentset{\circ} \sigma_\gamma = \overline \sigma_\gamma = \widehat \sigma_\gamma$. As for the drift, notice that since $\varphi$ is a chart in $M$, minimising $c_2$ will only involve a condition on the tangential part of $\accentset{\circ} \mu$, and is thus an unconstrained optimisation problem (the constraint \eqref{itotan} is then fulfilled by separately adding the required orthogonal term). Proceeding as in \autoref{sec:optimal}, we see that the quantity to be minimised is given by
\begin{equation}
	\sum_{p = 1}^m \bigg( \frac{\partial \overline\varphi^p}{\partial x^k} \accentset{\circ} \mu^k - \frac{\partial \overline\varphi^p}{\partial x^k} \frac{\partial \pi^k}{\partial x^h} \mu^h - \frac 12 \sum_{\gamma = 1}^n \frac{\partial^2 (\overline\varphi^p \circ \pi)}{\partial x^i \partial x^j} \sigma^i_\gamma \sigma^j_\gamma  \bigg)^2
\end{equation}
which results in 
\begin{equation}
	\frac{\partial \overline\varphi}{\partial x^k} \accentset{\circ} \mu^k = \frac{\partial \overline\varphi}{\partial x^k} \frac{\partial \pi^k}{\partial x^h} \mu^h + \frac 12 \sum_{\gamma = 1}^n \frac{\partial \overline\varphi}{\partial x^h} \frac{\partial^2 \pi^h}{\partial x^i \partial x^j} \sigma^i_\gamma \sigma^j_\gamma + \frac{\partial^2 \overline\varphi}{\partial x^i \partial x^j} \overline \sigma^i_\gamma \overline \sigma^j_\gamma
\end{equation}
Since the last term vanishes by \eqref{d2phi} we have that this formula for $\accentset{\circ} \mu$ coincides with the It\^o-jet projection $\widehat \mu$.
\begin{rem}[Optimality for Riemannian ambient manifolds]
	This reformulation of the optimality criterion allows us to generalise the statement of \autoref{maintheorem} to the case where $\bbR^d$ is substituted with a general Riemannian manifold $D$, of which $M$ is a Riemannian submanifold, \eqref{sde} with a diffusion-type SDE on $D$ (in any one of the three equivalent formulations), and the squared Euclidean norm in \eqref{projmin} is substituted with $\mathscr d_D(Y_t,\pi(X_t))^2$. By considering a Nash embedding of $D$ (and hence, transitively, of $M$) in $\bbR^r$ for large enough $r$, and extending the diffusion to a diffusion in $\bbR^r$, we have that the $\mathscr d_M$-optimal projection and the $\mathscr d_{\bbR^d}$-optimal projection both coincide with the It\^o-jet projection. But since $\mathscr d_{\bbR^d} \leq \mathscr d_D \leq \mathscr d_M$, this projection must also be $\mathscr d_D$-optimal, as is immediate by comparing Taylor expansions.
	
	We may also ask whether \autoref{maintheoremV} admits a generalisation to the Riemannian case. This can be done by substituting the difference $Y_t-X_t$ with $\psi(Y_t)-\psi(X_t)$ in both \eqref{obvmin} and \eqref{weakerror}, where $\psi$ is any normal chart for the ambient Riemannian manifold $D$ centred at the initial condition $y_0$, and the radius $r$ appearing in \eqref{tau} is chosen so that the ball of radius $r$ centred in $y_0$ is contained in the domain of $\psi$. The proof of optimality is straightforward from \autoref{maintheoremV} and the fact that $T_{y_0}\psi$ is a linear isometry, thus making the square
	\begin{equation}
		\begin{tikzcd}[column sep = large]
			T_{y_0}D \arrow[r,"T_{y_0}\psi"]\arrow[d,"P(y_0)"] & T_{y_0}\bbR^d \arrow[d,"P'(y_0)"] \\
			T_{y_0}M \arrow[r,"T_{y_0}\psi|_M"] & T_{y_0} M'
		\end{tikzcd}
	\end{equation}
	(where $M' = \mathrm{Im}\psi$, and $P(y_0),P'(y_0)$ are the metric projections) commute.
\end{rem}
We have thus shown that both \autoref{maintheoremV} and \autoref{maintheorem} can be reformulated so as to apply to the case of the ambient manifold being Riemannian.

\subsection{Optimality criteria for the Stratonovich projection}\label{subsec:strat}
It is surprising that the most na\"ive way to project the coefficients of an SDE is suboptimal according to the criteria introduced in this chapter. In this subsection we a (somewhat less compelling) way in which the Stratonovich projection can be considered optimal. This idea is already present in \cite[\S 4.4]{AB17}.

As before, we start with the Stratonovich SDE \eqref{StratSDE}. Define a second SDE
\begin{equation}\label{reflectedXi}
	\dif \Xi_t = -\sigma_\gamma(\Xi_t,t) \circ \dif B^\gamma_t - b(\Xi_t,t) \dif t, \quad \Xi_0 = y_0
\end{equation}
where $B$ is another $n$-dimensional Brownian motion, with no specific relationship with $W$. Assume we are looking for coefficients $\accentset{\circ} \sigma_\gamma$ and $\accentset{\circ} b$ s.t., defining
\begin{align}
	\begin{split}
		\dif Y_t &= \accentset{\circ} \sigma_\gamma(Y_t,t) \circ \dif W^\gamma_t + \accentset{\circ} b(Y_t,t) \dif t, \quad Y_0 = y_0 \\
		\dif \Upsilon_t &= -\accentset{\circ} \sigma_\gamma(\Upsilon_t,t) \circ \dif B^\gamma_t - \accentset{\circ} b(\Upsilon_t,t) \dif t, \quad \Upsilon_0 = y_0
	\end{split}	
\end{align}
the following quantity is minimised for small $t$ (in the same sense as in \autoref{maintheorem}):
\begin{equation}\label{symtoopt}
	\frac 12 E\big[|Y_t - \pi(X_t)|^2  \big] + \frac 12 E\big[|\Upsilon_t - \pi(\Xi_t)|^2  \big] 
\end{equation}
Note that the original input of the problem is the same as before, i.e.\ $\sigma_\gamma$ and $\mu$, but the quantity to be optimised is different. In \cite{ABR19} the SDEs with reflected Stratonovich coefficients are interpreted as representing a solution going backward in time: this fits in nicely with the interpretation of the Stratonovich integral of being time-symmetric (e.g.\ in the sense of the midpoint-evaluated Riemann sums that converge in $L^2$ to it). This interpretation is backed up by the fact that, if $\mu$ and $\accentset{\circ} \mu$ denote the It\^o drifts for $X$ and $Y$, the SDE for $\Xi, \Upsilon$ can be equivalently written using the backwards It\^o integral $\dif^\text{b}$ (defined by endpoint evaluation)
\begin{equation}
	\dif^\text{b} \Xi_t = -\sigma_\gamma(\Xi_t,t) \dif^\text{b} B^\gamma_t - \mu(\Xi_t,t) \dif t, \qquad \dif^\text{b} \Upsilon_t = -\accentset{\circ}\sigma_\gamma(\Upsilon_t,t) \dif^\text{b} B^\gamma_t - \accentset{\circ} \mu(\Upsilon_t,t) \dif t
\end{equation}
so that \eqref{symtoopt} can be viewed as averaging an SDE going forward in time with one going backwards. This heuristic interpretation, however, is not necessary in the computations, and we can proceed by optimising \eqref{symtoopt} as is. Proceeding as above, this leads to the the diffusion coefficients being, as always, orthogonally projected ($\accentset{\circ}\sigma_\gamma = \overline \sigma_\gamma$) and the constrained optimisation problem for the drift $\widetilde \mu$ given by
\begin{equation}
	\begin{dcases}
		\text{minimise }&\frac 12 \sum_{k = 1}^m \bigg(\widetilde \mu^k - \frac{\partial \pi^k}{\partial x^h}\mu^h - \frac 12 \sum_{\gamma = 1}^n \frac{\partial^2 \pi^k}{\partial x^i \partial x^j} \sigma^i_\gamma \sigma^j_\gamma \bigg)^2 \\
		& \!\!\!\!\!\!\!\!\!\!\!\!\!\!\!\!\!\!\!\!\!\!\!\!\!\!\!+ \frac 12 \sum_{k = 1}^m \bigg(\!\!-\widetilde \mu^k +\sum_{\gamma = 1}^n \overline \sigma^l_\gamma \frac{\partial \overline \sigma^k_\gamma}{\partial x^l} - \frac{\partial \pi^k}{\partial x^h}\bigg( \!\! - \mu^h + \sum_{\gamma = 1}^n \sigma^l_\gamma \frac{\partial \sigma^h_\gamma}{\partial x^l}  \bigg) - \frac 12 \sum_{\gamma = 1}^n \frac{\partial^2 \pi^k}{\partial x^i \partial x^j} \sigma^i_\gamma \sigma^j_\gamma \bigg)^2  \\
		\text{subject to } &Q^k_h \widetilde \mu^h - \frac 12 \sum_{\gamma = 1}^n \frac{\partial^2 \pi^k}{\partial x^i \partial x^j} \overline \sigma^i_\gamma \overline \sigma^j_\gamma = 0
	\end{dcases}
\end{equation}
which is checked, by using Lagrange multipliers as above, to have solution the Stratonovich-projected drift \eqref{stratdrift}. Therefore, the Stratonovich projection is optimal in this \say{time-symmetric} sense.

\starsection{Conclusions and further directions}\label{sec:projConcl}
\addcontentsline{toc}{section}{Conclusions and further directions}

In this chapter we have shown, re-expressing and improving on the ideas of \cite{ABR19}, how a concrete optimisation problem involving SDEs points towards the use of It\^o calculus on manifolds, while the results given by adopting the more commonly used Stratonovich calculus are suboptimal. 

It would be interesting to extend this optimisation result to the case where the equation is driven by $(1/4,1)\ni H$-fractional Brownian motion, in the sense of rough paths (which for $H > 1/2$ means in the sense of Young). Although this would amount to a generalisation of a Stratonovich equation, as seen in \eqref{btable} the optimal coefficients can still be expressed as a function of the original ones and their derivatives, and similar formulae could be shown to hold in the case of fractional noise. The rough-Skorokhod conversion formula \cite{CL19,CL20} could be of help here, although it would have to first be extended to cover the case in which the RDE has drift for the problem to be interesting.

\bibliographystyle{alpha} 
\renewcommand\bibname{\sc References}
\bibliography{merged}\addcontentsline{toc}{chapter}{References}

\end{document}